\documentstyle[11pt]{article}

\newtheorem{thm}{Theorem}[section]
\newtheorem{lem}[thm]{Lemma}
\newtheorem{cor}[thm]{Corollary}
\newtheorem{pro}[thm]{Proposition}

\newtheorem{defi}[thm]{Definition}

\newcommand{\gm }{\Gamma }
\newcommand{\lon }{\longrightarrow }
\newcommand{\be }{\begin{eqnarray*}}
\newcommand{\ee }{\end{eqnarray*}}

\newcommand{\per }{\backl }

\input amssym.def
\input amssym

\newcommand{\pf}{\noindent{\bf Proof.}\ }
\newcommand{\qed}{\begin{flushright} $\Box$\ \ \ \ \ \
                  \end{flushright}}
\newcommand{\complex}{{\Bbb C}}
\newcommand{\reals}{{\Bbb R}}

\newcommand{\frakg}{{\frak g}}
\newcommand{\frakh}{{\frak h}}

\newcommand{\half}{\frac{1}{2}}

\newcommand{\backl}{\mathbin{\vrule width1.5ex height.4pt\vrule
height1.5ex}}
\newcommand{\cala}{{\cal A}}

\newcommand{\cald}{{\cal D}}

\newcommand{\calm}{{\cal M}}

\newcommand{\smalcirc}{\mbox{\tiny{$\circ $}}}

\newcommand{\flb}{{[\![}}
\newcommand{\frb}{{]\!]}}
\newcommand{\eqq}{=}

\def\description label#1{\hfil\bf[#1]\hfil}
\newcommand{\parr}[1]{\frac{\partial  #1}{\partial \lambda^{i}}}

\newcommand{\rr}{\Lambda^{\#}}

\newcommand{\ih}{\frac{i}{\hbar}}

\newcommand{\Hoch}{Hochschild }

\newcommand{\ot}{\mbox{$\otimes$}}

\newcommand{\otr}{\ot_{R}}

\newcommand{\oh}{O(\hbar^{2})}

\newcommand{\alt}{\mbox{Alt}}

\newcommand{\vh}[1]{\Vec{h_{#1}}}
\newcommand{\ve}[1]{\Vec{e_{#1}}}
\newcommand{\lam}[1]{\frac{\partial}{\partial \lambda^{#1}}}

\newcommand{\third}{\frac{1}{3}}
\newcommand{\muu}{{\frak m}}
\newcommand{\nablaa}{\nabla^{0}}
\newcommand{\hs}[1]{h_{*}^{#1}}
\newcommand{\es}[1]{e_{*}^{#1}}

\newcommand{\vhs}[1]{\Vec{h_{*}^{#1}}}
\newcommand{\ves}[1]{\Vec{e_{*}^{#1}}}

\newcommand{\ugh}{U\frakg \flb\hbar \frb}
\newcommand{\ug}{U\frakg}
\newcommand{\vT}{\Vec{\Theta}}
\newcommand{\vFT}{\Vec{F(\lambda )\Theta }}
\newcommand{\vF}{\Vec{F(\lambda )}}
\newcommand{\vdiT}{\Vec{[(\Delta \ot id )\Theta ] \Theta^{12} }}
\newcommand{\vidT}{\Vec{[(id \ot  \Delta) \Theta ]\Theta^{23}  }}
\newcommand{\perr}{\perp}
\newcommand{\calt}{{\cal T}}

\newcommand{\r}{\gamma}
\newcommand{\vs}[1]{v_{*}^{#1}}
\newcommand{\us}[1]{u_{*}^{#1}}
\newcommand{\A}{\cala}
\newcommand{\DD}{{\cal D}}

\def\sdp{\mathbin{\hbox{$\mapstochar\kern-.3333em\times$}}}
\def\pds{\mathbin{\hbox{$\times\kern-.55em\mapstochar\,$}}}

\newcommand{\wed}{\mathbin{\lower1.5pt\hbox{$\scriptstyle{\wedge}$}}}

\let\Tilde=\widetilde

\let\Vec=\overrightarrow

\def\chigh{{\raise1.5pt\hbox{$\chi$}}}
\let\phi=\varphi
\def\til0{\Tilde{0}}

\def\dminus{\raise2pt\hbox{\vrule height1pt width 2ex}\hskip3pt}

\def\pback#1{\mathbin{{{\lower1.2ex\hbox{$\times$}}\atop #1}}}

\def\vlra{\hbox{$\,-\!\!\!-\!\!\!-\!\!\!-\!\!\!-\!\!\!
-\!\!\!-\!\!\!-\!\!\!-\!\!\!-\!\!\!\longrightarrow\,$}}

\def\gpd{\,\lower1pt\hbox{$\longrightarrow$}\hskip-.24in\raise2pt
             \hbox{$\longrightarrow$}\,}

\def\lgpd{\,\lower1pt\hbox{$\vlra$}\hskip-1.02in\raise2pt\hbox{$\vlra$}\,}

\def\llgpd{\,\lower1pt\hbox{$\vvlra$}\hskip-1.3in\raise2pt\hbox{$\vvlra$}\,
}

\begin{document}

\title{{\bf Triangular dynamical $r$-matrices and quantization}}
\author{ PING XU \thanks{ Research partially supported by NSF
        grants   DMS97-04391 and DMS00-72171.}\\
 Department of Mathematics\\The  Pennsylvania State University \\
University Park, PA 16802, USA\\
        {\sf email: ping@math.psu.edu }}

\date{}

\maketitle
\begin{abstract}
We   study some general aspects of   triangular dynamical
$r$-matrices using Poisson geometry.  We show that a triangular dynamical
$r$-matrix  $r: \frakh^* \lon \wedge^{2}\frakg$ always  gives rise to a regular Poisson manifold.  Using the 
 Fedosov method, we prove that non-degenerate 
triangular dynamical $r$-matrices (i.e., those  such that the
 corresponding Poisson manifolds are symplectic) are quantizable, and that
the quantization is classified by  the     relative  Lie algebra cohomology
$H^{2}(\frakg , \frakh )\flb\hbar \frb$.
\end{abstract}

\section{Introduction}

In the last two decades, the theory of quantum groups
has undergone  tremendous development.  The classical
counterparts of quantum groups are Lie bialgebras \cite{Dr2}.
Many interesting quantum groups were found and studied
by various authors,  but
the proof of existence    of quantization for arbitrary
Lie bialgebras was obtained only  recently by
Etingof and  Kazhdan \cite{EK}. For triangular
Lie bialgebras, however, an  elementary
proof of quantization was given by Drinfeld
in 1983 \cite{Dr1}. Drinfeld's idea can be outlined
as follows. A triangular $r$-matrix on a Lie algebra
$\frakg$ defines a left invariant  Poisson structure on its corresponding
Lie group $G$. By restricting  to a Lie subalgebra if
necessary, one may in fact assume that this is 
symplectic. One may  then quantize
the $r$-matrix by finding a  $G$-invariant
$*$-product  on $G$, of which there may be several.
In \cite{Dr1},  Drinfeld   identified 
the symplectic manifold with a coadjoint orbit
of a  central extension of $\frakg$, and  then  applied
Berezin quantization \cite{B}.

Recently, there  has been 
  growing interest in the so-called quantum dynamical
Yang-Baxter equation (see Equation (\ref{eq:qdybe})).
This equation  arises
naturally from  various contexts in mathematical physics.
It first appeared in the work of Gervais-Neveu in their
study of quantum Liouville theory \cite{GN}.  Recently it reappeared
in Felder's work on the  quantum Knizhnik-Zamolodchikov-Bernard
 equation. It also  has  been found
to be connected  with the quantum Caloger-Moser systems
 \cite{ABE}.
Just like  the  quantum Yang-Baxter equation is connected
with quantum groups, the  quantum dynamical
Yang-Baxter equation is known to be connected with elliptic
quantum groups \cite{Felder}, 
as well as with Hopf algebroids or quantum groupoids \cite{EV2, EV3, Xu2, Xu:gpoid}.

 The classical counterpart
of the quantum dynamical Yang-Baxter equation  was first
considered by Felder \cite{Felder}, and then studied
by Etingof and Varchenko \cite{EV1}. This is 
 the so-called classical dynamical
Yang-Baxter equation, and a solution to such an
equation (plus some other reasonable conditions) is
called a classical dynamical $r$-matrix. More precisely,
given a Lie algebra $\frakg$ over $\reals$ (or over $\complex$)
 with an Abelian Lie subalgebra $\frakh$, a classical dynamical $r$-matrix 
is a smooth (or meromorphic)  function $r(\lambda ): \frakh^* \lon 
\frakg \ot \frakg $ 
satisfying the following conditions:
\begin{enumerate}
\item (zero weight condition) $[h\ot 1 +1 \ot h, r(\lambda )]=0, \  \ \forall h\in \frakh$;
\item (normal condition)  $r^{12}+r^{21}= \Omega$, where
$\Omega\in (S^{2}\frakg)^{\frakg}$ is a Casimir element;
\item (classical dynamical Yang-Baxter equation)
\begin{equation}
\label{cdybe}
 \alt (dr) \, + \,
[r^{12},r^{13}] + [r^{12}, r^{23}] + [r^{13}, r^{23}] \, = \, 0,
\end{equation}
\end{enumerate}
where $\alt dr =\sum (h_{i}^{(1)} \frac{\partial r^{23}}{\partial \lambda^{i}}
-h_{i}^{(2)} \frac{\partial r^{13}}{\partial \lambda^{i}}+
h_{i}^{(3)} \frac{\partial r^{12}}{\partial \lambda^{i}} )$.

 A fundamental question is whether any  classical dynamical $r$-matrix
is  quantizable.  There have appeared many
results in this direction. For  the standard classical dynamical $r$-matrix
for  $\frak{sl}_{2}(\complex )$, a  quantization was
  obtained by Babelon  \cite{Babelon} in 1991.
For general simple Lie algebras, quantizations were recently
 found  independently by Arnaudon et al. \cite{Arnaudon}  and Jimbo
et al. \cite{Jimbo}  based on the approach of Fronsdal \cite{F}.
Similar results  were also found by Etingof and
Varchenko \cite{EV3} using intertwining operators.
Recently, using a  method similar to \cite{Arnaudon, F, Jimbo},
 Etingof et al. \cite{ESS} obtained a quantization 
of  all  the classical dynamical $r$-matrices of semi-simple Lie algebras
in  Schiffmann's classification  list \cite{S}. However,
the general quantization problem  still remains  open;
 a   recipe has   yet  to be found.  Moreover, the problem
of classification of quantizations has not yet been touched.

In this paper, we study the  quantization problem for
general  classical triangular  dynamical $r$-matrices.
Classical triangular  dynamical $r$-matrices are those
satisfying the skew-symmetric condition
 $r^{12}(\lambda )+ r^{21}(\lambda )=0$.  In this case, Equation (\ref{cdybe}) 
is equivalent to  
 $\sum_{i} h_{i} \wedge \frac{\partial r}{ \partial\lambda^{i}} +\half [r , r]= 0$.
These $r$-matrices are  in one-one correspondence with regular
 Poisson structures  $\pi =\sum_{i} \vh{i}\wedge \parr{} +\Vec{r(\lambda )}$ 
on the manifold $\frakh^* \times G$, which are
 invariant under the left $G$ and
right $H$-actions.
Thus one may expect to quantize a classical dynamical $r$-matrix
by looking for 
 a certain special type of  star-products \cite{BEFFL} on the
 corresponding Poisson manifold. This is
exactly the route we take in the present paper. In some
sense, this is also a natural generalization of the
quantization method used by Drinfeld in \cite{Dr1} as outlined
at the beginning  of the introduction. In fact, in the
present paper, we mainly deal with non-degenerate triangular
 classical dynamical $r$-matrices (i.e., the corresponding
 Poisson manifolds are in fact symplectic).  Berezin
quantization no longer works in this  situation.  However,
 one may use the Fedosov method to obtain the desired star-products
as we will see later.
It is well-known that star products on a symplectic
manifold are classified by the  second cohomology
group of the manifold with coefficients in  formal
$\hbar$-power series. In light of this result, we are
able to classify the quantizations of a non-degenerate triangular
 classical dynamical $r$-matrix and prove that the
 quantizations are  parameterized by the relative Lie algebra
cohomology $H^{2}(\frakg ,\frakh )\flb\hbar \frb$.

For a general  triangular classical dynamical $r$-matrix, it is natural
to  ask whether it is  possible to reduce it to
a non-degenerate one by restricting  to   a Lie subalgebra.
 This is always true in the non-dynamical case
\cite{Dr1}. Unfortunately, in general this fails in the dynamical case,
and we will study the conditions  under which this is possible. In this case,
these
$r$-matrices are called splittable.
Splittable triangular classical dynamical $r$-matrices resemble in many
ways  non-degenerate ones. And in particular, they can be
quantized by  the Fedosov method.

The outline of this paper is  as follows. After Section 1 (this introduction),
  in Section 2, we study  general
properties of triangular classical dynamical $r$-matrices.
It is proved that triangular classical dynamical $r$-matrices
correspond to  some special Poisson structures on $\frakh^* \times G$,
which are always regular. This may seem surprising at first
glance since the rank of $r(\lambda )$ may depend on the point
$\lambda$. The main tool in Section 2
is the method of Lie groupoids and Lie algebroids. In particular, we show how  gauge transformations,
  first introduced by Etingof and Varchenko  \cite{EV1},
enter naturally from the viewpoint of  Lie algebroids.  
The study of the tangent space of the moduli  space of
dynamical $r$-matrices naturally
leads  to  the notion of
 dynamical $r$-matrix cohomology, which is shown
to be isomorphic to the relative Lie algebra cohomology when
 $r$ is non-degenerate. Section 3 is devoted to the proof  of the
equivalence between quantizations of triangular classical dynamical $r$-matrices
and the so called compatible star products on their
corresponding Poisson manifolds $\frakh^* \times G$. In Section 4, we study
 symplectic connections on such  symplectic manifolds ($M=\frakh^* \times G$).
In particular, we show that there always exists a $G\times H$-invariant
(i.e. left $G$-invariant and right $H$-invariant) torsion-free
symplectic  connection on $M$  such that the left invariant vector
fields $\vh{}, \ \forall h\in \frakh$ are all parallel.
The main result of Section 5 is 
that the Fedosov quantization obtained via such   a symplectic connection
and some suitable choice of Weyl curvatures gives rise to compatible
$*$-products on $M=\frakh^* \times G$. Therefore, as a 
consequence, we prove  the existence of a  quantization 
of non-degenerate triangular classical dynamical $r$-matrices.
The presentation in Section 5, however, is made in a more general
setting, which is  of its own interest.
Section 6 is devoted to the classification of quantizations.
In particular, we show that  the
 equivalence classes of quantizations of a non-degenerate triangular
 dynamical $r$-matrix
$r:\frakh^* \lon \wedge^2  \frakg$ are parameterized by
the relative Lie algebra cohomology  with 
coefficients in the formal $\hbar$-power series
$H^{2}(\frakg , \frakh )\flb\hbar \frb$.  Some speculation
on the classification of quantizations of a general
triangular classical dynamical $r$-matrix is
given as a conjecture, which is consistent with  Kontesvich's formality
theorem \cite{K}.
In the appendix we  recall some basic ingredients  of
the Fedosov  quantization, which are used throughout the paper.

Finally, some  remarks are in order. Quantization of
dynamical $r$-matrices is related to quantization
of Lie bialgebroids as shown in \cite{Xu:gpoid}.
However, for simplicity, we will avoid using  quantum groupoids
in the present paper even though many ideas  are rooted from
  there.
   Also in this paper, we work in the smooth case. Namely,
 Lie algebras are  finite dimensional  Lie algebras  over $\reals$,
all manifolds and maps are  smooth,
 but our approach works for the
 complex category as well.
For simplicity,  we assume
that a  dynamical $r$-matrix is always  defined on $\frakh^*$. In reality,
it may only be defined on an open submanifold
  $U\subset \frakh^*$,  but  our results   hold in this situation as well.

{\bf Acknowledgments.} The author  would   like to thank  
  Martin Bordemann,  Pavel Etingof and Boris Tysgan
 for useful discussions.
Especially, he is  grateful to  Pavel Etingof for
his  suggestion of writing up  this work.
 In addition to the funding sources mentioned
in the first footnote, he would also  like to thank
  the Max-Planck Institut for the
 hospitality and financial support while part of this project was being done.

\section{Triangular dynamical $r$-matrices}
In this section, we study some  general aspects of
triangular dynamical $r$-matrices.
As a useful tool, we  shall  utilize
 the method of Lie algebroids and
Lie groupoids. Let $\frakg$ be a Lie algebra and $\frakh\subset \frakg$
an Abelian Lie subalgebra of dimension $l$.
 By a {\it  triangular dynamical $r$-matrix}, we mean
a  smooth  function  $r:\frakh^{*} \lon \wedge^2 \frakg  $
satisfying:
\begin{enumerate}
\item the zero weight condition:
 $[h, r (\lambda )] = 0, \ \ \ \forall \lambda\in \frakh^*, \ h \in \frakh $, and
\item   the classical dynamical Yang-Baxter equation (CDYBE):
\begin{equation}
\label{eq:cdybe}
\sum_{i} h_{i} \wedge \frac{\partial r}{ \partial\lambda^{i}} +\half [r , r]=0,
\end{equation}
\end{enumerate}
where the bracket $[\cdot , \cdot ]$ refers to the Schouten type
bracket: $\wedge^{k}\frakg \ot \wedge^{l}\frakg \lon \wedge^{k+l-1}\frakg $
induced from the Lie algebra bracket on $\frakg$.
Here $\{h_{1}, \cdots , h_{l}\}$  is a basis in  $\frakh$,  and
$(\lambda^{1}, \cdots , \lambda^{l})$ its  induced
coordinate system  on $\frakh^*$. It is known \cite{BK-S} \cite{LX1}
that the CDYBE is closely related to Lie bialgebroids.
 Recall that  a Lie bialgebroid is a pair of Lie algebroids 
($A$, $A^*$) satisfying
 the following compatibility condition (see \cite{MX94, MX98, K-S95}):
\begin{equation}
\label{eq:1}
d_{*}[X, Y]=[d_{*}X, Y]+[X, d_{*}Y], \ \ \forall X ,Y \in \Gamma (A),
\end{equation}
where the differential $d_*$ on $\Gamma(\wedge^*A)$
comes from the Lie algebroid structure on $A^*$.

Given  a Lie algebroid $A$  over $P$ with anchor $a$, and a section
$\Lambda $ of $ \Gamma (\wedge^2 A )$ satisfying
the condition $[\Lambda, \Lambda ]=0$, one may
define a Lie algebroid structure on $A^*$ by
simply requiring the differential $d_* : \Gamma(\wedge^k A)\lon
\Gamma(\wedge^{k+1} A)$ to be $d_*  =[\Lambda , \cdot ]$.
 More explicitly, denote by $\Lambda^\#$
 the bundle map $A^{*}\lon A$ defined by
$\rr (\xi) (\eta )=\Lambda  (\xi , \eta ), \forall \xi , \eta \in
 \Gamma (A^{*})$.  Then the  bracket on
$\Gamma (A^* )$  is defined by
\begin{equation}
\label{eq*}
[\xi , \eta ] = L_{\rr \xi}\eta -
L_{\rr \eta}\xi -d[\Lambda  (\xi , \eta )],
\end{equation}
and the anchor $a_*$ is the composition  $a\circ \rr :A^*  \lon TP$.
It is easy to show that $(A, A^* )$ is  indeed a 
Lie bialgebroid, which is called a {\em triangular
Lie bialgebroid} \cite{MX94}.

Now consider $A=T\frakh^* \times \frakg$ and equip $A$ with
the standard  product Lie algebroid structure.
Then the anchor $a:T\frakh^* \times \frakg \lon T\frakh^* $
is  simply the projection. The relation between
triangular dynamical $r$-matrices and
triangular Lie bialgebroids are described by the following 
\cite{BK-S, LX1}:

\begin{pro}
\label{pro:triangular}
Given a smooth function  $r:\frakh^{*} \lon \wedge^2 \frakg  $,
$r$ is a triangular dynamical $r$-matrix iff the Lie algebroid
$(A, a)$ together with $\Lambda =\sum_{i} h_{i} \wedge \parr{}
 +r(\lambda )\in \gm ( \wedge^2  A)$ defines a triangular Lie bialgebroid.
\end{pro}
\pf  By   a straightforward computation, we have  $[\Lambda , \Lambda ]
=2(\sum_{i}  h_{i} \wedge \frac{\partial r}{ \partial\lambda^{i}} +\half [r, r]
+\sum_{i} [r, h_{i}]\wedge \parr{} )$. It thus 
 follows that $[\Lambda , \Lambda ]=0$ iff
 $\sum_{i}  h_{i} \wedge \frac{\partial r}{ \partial\lambda^{i}}
 +\half [r, r]=0$ and $[r, h_{i}]=0 \ (i=1, \cdots l$), i.e.,
 $r$ is a triangular dynamical 
$r$-matrix. \qed

Let $G$ be a Lie group with  Lie algebra $\frakg $ and $H\subset G$ an Abelian Lie subgroup
with Lie algebra $\frakh$.
 Consider  $M=\frakh^* \times G$. Let $G$ act on $M$  from the left by left multiplication
on $G$, and  $H$ act from the right by  right multiplication on $G$.
An equivalent version of Proposition \ref{pro:triangular}
is the following
\begin{pro}
\label{pro:Poisson}
For a smooth function  $r:\frakh^{*} \lon \wedge^2 \frakg  $,
$r$ is a triangular dynamical $r$-matrix iff 
$\pi =\sum_{i} \vh{i}\wedge \parr{} +\Vec{r(\lambda )}$ defines a $G\times H$-invariant Poisson
structure on $M=\frakh^* \times G$, where $\vh{i}\in {\frak X} (M) $ is the left invariant
vector field on $M$ generated by $h_{i}$ and  similarly
$\Vec{r(\lambda )}\in \gm (\wedge^{2} TM)$ is the left invariant bivector field
on $M$ corresponding to  $r(\lambda )$.
\end{pro}

\begin{thm}
\label{thm:iso}
If $r:\frakh^{*} \lon \wedge^2 \frakg  $ is a triangular dynamical $r$-matrix,
then $\frakh +r(\lambda )^{\#}\frakh^{\perp}$ is a 
Lie subalgebra of $\frakg$.  Moreover 
the Lie subalgebras $\frakh +r(\lambda )^{\#}\frakh^{\perp}$, $\forall
\lambda \in \frakh^{*}$, 
  are all  isomorphic,  and the isomorphisms are given by
  the adjoint action of $G$.
\end{thm}
\pf  For any $\lambda \in \frakh^{*}$, $A_{\lambda}=T_{\lambda}\frakh^{*}
\oplus \frakg\cong \frakh^{*} \oplus \frakg$ and $A_{\lambda}^*
\cong \frakh \oplus \frakg^*$. Under these identifications,
 the bundle map   $\Lambda^{\#}_{\lambda}:
A_{\lambda}^* \lon A_{\lambda}$ is given by
\begin{equation}
\label{eq:bundle}
 (h, \xi )\mapsto (i^* \xi , -h+r(\lambda )^{\#} \xi ),  \ \ \forall h\in \frakh \ \mbox{  and }
\xi \in \frakg^*,
\end{equation}
where $i:\frakh \lon \frakg$ is the inclusion.
 Set $B=\Lambda^{\#} (A^* )=\cup_{\lambda\in \frakh^* }
\Lambda^{\#}_{\lambda} (A_{\lambda}^* ) \subset A$.
 Since $(A, \ \Lambda)$ defines a triangular Lie bialgebroid,
$B$ is integrable. I.e., $\gm (B)$ is closed under the 
Lie algebroid bracket on $\gm (A)$. Hence $\mbox{ker} a|_{B_{\lambda}}$
is  a Lie subalgebra of $\mbox{ker} a|_{A_{\lambda}}$. Now
it is easy to see that $\mbox{ker} a|_{B_{\lambda}}=\frakh +r(\lambda )^{\#}\frakh^{\perp}$
 and $\mbox{ker} a|_{A_{\lambda}}=\frakg$.  It thus follows that
$\frakh +r(\lambda )^{\#}\frakh^{\perp}$ is a
Lie subalgebra of $\frakg$. On the other hand,  from Equation (\ref{eq:bundle}), it
is easy to see that $a  (B_{\lambda} ) =T_{\lambda} \frakh^*$. Hence
$a: B\lon T\frakh^{*} $ is surjective, which implies that
$B$ is  in fact a transitive Lie
algebroid (also called a gauge Lie algebroid \cite{Mackenzie}). Thus
it follows that  the dimension  of $B_{\lambda} $ is  independent
of $\lambda$,  and therefore  $B$ is a subbundle of $A$.
Moreover  the isotropic Lie algebras of $B$ at different
points of $\frakh^*$ are all isomorphic, 
and the isomorphisms are given by the
adjoint action of $G$.
This implies that, for any $\lambda , \mu \in \frakh^{*}$, 
$\frakh +r(\lambda )^{\#}\frakh^{\perp}$   is isomorphic to
$\frakh +r(\mu )^{\#}\frakh^{\perp}$  by the  adjoint action of
a group element in $G$.  \qed

For the sake of   simplicity,
 we denote by $\frakg_{\lambda}$
the Lie subalgebra $\frakh +r(\lambda )^{\#}\frakh^{\perp}$.
Define the {\em rank of a triangular dynamical $r$-matrix} $r $  to be
$\mbox{dim} \frakg_{\lambda}- \mbox{dim} \frakh$, which is denoted as
rank$r$.
We say a triangular dynamical $r$-matrix  $r$ is
 {\em non-degenerate} if rank$r=\mbox{dim} \frakg- \mbox{dim} \frakh$.

An immediate consequence of Theorem \ref{thm:iso} is

\begin{cor} 
Under the same hypothesis as in Theorem \ref{thm:iso}, rank$r$
is independent of the point $\lambda$ and  therefore is  a
well-defined even number. Moreover
$B=\Lambda^{\#}A^{*} \subset A$ 
is a Lie subalgebroid of  rank $2\mbox{dim}\frakh $+rank$r$, and
$(M, \pi )$ is a regular Poisson manifold of rank  
$2\mbox{dim}\frakh $+rank$r$.
\end{cor}

In particular, we have the following
\begin{cor}
\label{cor:sym}
Given a triangular dynamical $r$-matrix $r: \frakh^* \lon \wedge^2 \frakg$,
 the following  statements are all equivalent:
\begin{enumerate}
\item $r$ is   non-degenerate;
\item the bundle map  $\Lambda^{\#}: A^* \lon A$ is nondegenerate;
\item $\frakg_{\lambda}=\frakg, \ \ \forall \lambda \in \frakh^*$;
\item $(M, \pi)$ is a symplectic manifold.
\end{enumerate}
\end{cor}

If we  choose a decomposition $\frakg =\frakh \oplus \muu$,
where $\muu$ is a    subspace of  $ \frakg$, 
 and choose a basis $\{h_{1}, \cdots , h_{l} \}$ for
$\frakh$ and a basis $\{e_{1},  \cdots , e_{m}\}$ for $\muu$,
we may write
\begin{equation}
\label{eq:dec}
r(\lambda )=\sum a^{ij}(\lambda )h_{i}\wedge h_{j} +
\sum b^{ij}(\lambda )h_{i}\wedge e_{j} + \sum c^{ij}(\lambda )e_{i}\wedge e_{j}.
\end{equation}
It is simple to see that $\frakg_{\lambda}=\frakh \oplus \mbox{Span}\{
\sum_{j}c^{ij}(\lambda )e_{j}|i=1, \cdots , m\}$, and rank$r$ is
the rank of the matrix $(c^{ij}(\lambda ))$. Therefore,
we immediately know  that the rank of $(c^{ij}(\lambda ))$
is independent of $\lambda $. Clearly  $r$ is
non-degenerate iff the matrix $(c^{ij}(\lambda ))$ is
non-degenerate.

A  natural question arises as to whether it is possible
to make an arbitrary triangular dynamical $r$-matrix non-degenerate
by considering it to be valued in a Lie subalgebra of
$\frakg$.
 This is true in the non-dynamical case \cite{Dr1}, for example.
However, in the dynamical case, this is not always 
possible as we will see below. Nevertheless we will
single out those $r$-matrices possessing this property, which
will be called {\em splittable}. 
 Splittable triangular dynamical $r$-matrices
contain a large class of  interesting dynamical $r$-matrices, which
in fact include almost  all  examples we know,
e.g., those as classified in \cite{EV1} when $\frakg$ is a simple
Lie algebra. More precisely,
\begin{defi}
A triangular dynamical $r$-matrix $r: \frakh^* \lon \wedge^2 \frakg$
 is said to be splittable
if for any $\lambda \in \frakh^*$, $i^* (r(\lambda )^{\#-1} \frakh )=
\frakh^*$, where $i: \frakh \lon \frakg$ is the inclusion.
\end{defi}

\begin{pro}
\label{pro:split}
Suppose that $r$ is a  triangular dynamical $r$-matrix. Then the  following
statements  are equivalent:
\begin{enumerate}
\item $r$ is  splittable;
\item for any $\lambda \in \frakh^*$,
 $ r(\lambda )^{\#} \frakg^* \subset \frakg_{\lambda}$;
\item if $r(\lambda )$ is given as in Equation (\ref{eq:dec}) under
 a decomposition $\frakg =\frakh \oplus \muu$,
then for any $i$, $\sum_{j}b^{ij}(\lambda )e_{j}\in 
\mbox{Span}\{
\sum_{j}c^{ij}(\lambda )e_{j}|i=1, \cdots , m\}$;
\item for any fixed $\lambda \in \frakh^*$, there exists a
decomposition $\frakg =\frakh \oplus \muu$, under which
\begin{equation}
\label{eq:uniform}
r(\lambda )=\sum a^{ij}(\lambda )h_{i}\wedge h_{j} + \sum c^{ij}(\lambda )e_{i}\wedge e_{j};
\end{equation}
\item $T\frakh^* \times \{0\}\subset B$.
\end{enumerate}
\end{pro}

Let us first prove the following simple lemma  from 
linear algebra. 

\begin{lem}
\label{lem:basis-change}
Let $V=\frakh \oplus \muu$ be a decomposition of vector spaces,
and let $\{ h_{1}, \cdots , h_{l}\}$ be a basis of $\frakh$, and
$\{e_{1},  \cdots , e_{m}\}$ a basis of   $\muu$.
Let $r\in \wedge^2 V$ be any element such that
 $$r=\sum a^{ij} h_{i}\wedge h_{j} +
\sum h_{i}\wedge x^{i} + \sum c^{ij}e_{i}\wedge e_{j}, $$
where $x^{i}\in \muu$, and  $a_{ij}, \ c_{ij}$ are
skew-symmetric, i.e., $a_{ij}=-a_{ji}$ and
$c_{ij}=-c_{ji}$. If $I\subset \{1, \cdots , l\}$
 is a subset of indexes such that
for any $i_{0}\in I$,  $x^{i_{0}}\in {Span}\{
\sum_{j}c^{ij}e_{j}|i=1, \cdots , m\}$. Then one can change 
the decomposition 
$V=\frakh \oplus \tilde{\muu}$ so that under
a suitable basis  $\{\tilde{e_{1}},  \cdots , \tilde{e_{m}}\} $
of $\tilde{\muu}$,  $r$ can be written as
$$r=\sum \tilde{a}^{ij} h_{i}\wedge h_{j} +
\sum_{i\notin I} h_{i}\wedge x^{i} + \sum c^{ij}\tilde{e}_{i}\wedge 
\tilde{e}_{j}.  $$
\end{lem}
\pf $\forall i_{0}\in I$, by assumption,  there are constants 
$\r_{i}^{i_{0}}, \ i=1, \cdots , m$,
such that $x^{i_{0}} =2\sum_{ij} \r_{i}^{i_{0}} c^{ij}e_{j}$.
Let $ \tilde{e}_{i}=e_{i}+\sum_{i_{0}\in I }\r_{i}^{i_{0}}h_{i_{0}},
\ \ \forall i=1, \cdots , m$. Then
\be
&&\sum c^{ij}\tilde{e}_{i}\wedge \tilde{e}_{j}\\
&=&\sum c^{ij}(e_{i}+\sum_{i_{0}\in I }\r_{i}^{i_{0}}h_{i_{0}})\wedge
(e_{j}+\sum_{i_{0}\in I }\r_{j}^{i_{0}}h_{i_{0}})\\
&=&\sum c^{ij} e_{i}\wedge e_{j}+2 \sum c^{ij} \r_{i}^{i_{0}} h_{i_{0}} \wedge
e_{j} \ \ \ \ \ (\mbox{mod} \wedge^2 \frakh )\\
&=&\sum c^{ij}e_{i}\wedge e_{j}+\sum_{i_{0}\in I } h_{i_{0}}
\wedge x^{i_{0}} \ \ \ \ \ (\mbox{mod} \wedge^2 \frakh ).
\ee
Hence $r=\sum c^{ij}\tilde{e}_{i}\wedge \tilde{e}_{j}+\sum_{i\notin I }
h_{i} \wedge x^{i} \ \ \ (\mbox{mod} \wedge^2 \frakh )$.
This concludes the proof. \qed 
{\bf Proof of Proposition \ref{pro:split}} 

(i)$\Rightarrow$(ii) Let us fix a basis $\{h_{1}, \cdots , h_{l}\}$ of $\frakh$, and
let $\{h^1_*, \cdots , h^l_* \}$ be its dual basis in $\frakh^*$.
By assumption, for any $1\leq j\leq l$, there is a $\xi^j \in \frakg^*$
such that $i^* \xi^j =h^j_*$ and $r(\lambda )^{\#}\xi^j \in \frakh$.
Given any $\xi \in \frakg^*$, take $a_{j}=<\xi , h_{j}>$ and
$\eta =\xi -\sum a_{j}\xi^j$. Then it is easy to see that $\eta  \in 
\frakh^{\perp}$. Hence $r(\lambda )^{\#}\xi =
\sum a_{j} r(\lambda )^{\#}\xi^j
+ r(\lambda )^{\#}\eta \in \frakh +r(\lambda )^{\#}\frakh^{\perp} =
\frakg_{\lambda }$.

(ii)$\Rightarrow$(iii)  Let $\{ h_{1},  \cdots , h_{l}\}$ be
a basis of $\frakh $, $\{e_{1},
\cdots ,  e_{m} \}$ a basis of $\muu$,  and
$\{ h^{1}_{*}, \cdots , h^{l}_{*}, e^{1}_{*},
\cdots ,  e^{m}_{*} \}$  the dual basis of $\{ h_{1},  \cdots , h_{l}, e_{1},
\cdots ,  e_{m} \}$ in $\frakg^*$.
It is trivial to see that $r(\lambda )^{\#} e^{i}_{*}=-\sum_{j}
b^{ji}(\lambda )h_{j}+2 \sum_{j}c^{ij}(\lambda )e_{j}$. Hence we have
$$\frakg_{\lambda }=\frakh \oplus \mbox{Span}\{ \sum_{j}c^{ij}(\lambda )e_{j}|i=1, \cdots , m\}.$$
 Now $r(\lambda )^{\#}h^{i}_{*} =\sum_{j}2a^{ij}(\lambda )h_{j}+
\sum_{j} b^{ij}(\lambda )e_{j}$. Since $r(\lambda )^{\#}h^{i}_{*}\in
\frakg_{\lambda }$ by assumption, it follows  that $\sum_j b^{ij}(\lambda )e_{j}\in 
\mbox{Span}\{
\sum_{j}c^{ij}(\lambda )e_{j}|i=1, \cdots , m\}$.

(iii)$\Rightarrow$(vi)  This follows from Lemma  \ref{lem:basis-change}.

(vi)$\Rightarrow$(v) If $r(\lambda )
=\sum a^{ij}(\lambda )h_{i}\wedge h_{j} + \sum c^{ij}(\lambda )e_{i} 
\wedge e_{j}$, then $r(\lambda )^{\#}h^{i}_{*} =2\sum_{j} a^{ij}(\lambda )h_{j}$. Thus according to 
Equation (\ref{eq:bundle}), $\Lambda_{\lambda}^{\#} (2\sum_{j} 
a^{ij}(\lambda )h_{j},  \ h^i_* )=(h^i_* , 0)$. Hence,
  $(h^i_* , 0)\in B_{\lambda}$. This implies that $T_{\lambda}\frakh^* \times \{0\}\subset B_{\lambda}$.

(v)$\Rightarrow$(i) Given any  $\phi \in \frakh^*$,
we know that  $(\phi ,0)\in B_{\lambda}$ by assumption.
Therefore there exist  $h\in \frakh$ and $\xi \in \frakg^*$ such that
$\Lambda_{\lambda}^{\#} (h, \xi )=(\phi , 0)$, i.e., 
$(i^{*}\xi , -h +r(\lambda )^{\#}\xi )=(\phi, 0)$ according
to Equation (\ref{eq:bundle}). This
implies that $\phi =i^{*}\xi $ and $r(\lambda )^{\#}\xi =h$.
Hence $\phi \in  i^* (r(\lambda )^{\#-1} \frakh )$. Therefore, we conclude
that $  \frakh^* \subset  i^* (r(\lambda )^{\#-1} \frakh )$. \qed
{\bf Remark} In the  proof above, the decomposition
  $\frakg =\frakh \oplus \muu$ and
the choice of the basis $\{e_{1}, \cdots , e_{m}\}$
in (iv) depend on a   particular point $\lambda$. It is  not clear 
whether it is possible to find a decomposition so that 
Equation (\ref{eq:uniform}) holds uniformly for all points in $\frakh^*$.
 On the other hand, if there exists such  a decomposition 
$\frakg =\frakh \oplus \muu$ so that  a triangular dynamical $r$-matrix
is of the form as in Equation (\ref{eq:uniform}), it is always splittable. \\\\

An immediate consequence of Proposition \ref{pro:split} is the following:
\begin{cor}
\label{cor:split-nondeg}
If  $r:\frakh^*  \lon \wedge^{2}\frakg $ is a splittable
triangular  dynamical $r$-matrix, then
\begin{enumerate}
\item $\frakg_{\lambda}$ is independent of $\lambda$, i.e.,
$\frakg_{\lambda}=\frakg_{\mu}, \ \forall \lambda , \mu \in \frakh^*$.
We will denote $\frakg_{\lambda}$ by $\frakg_1$.
\item $r$ can be considered as a non-degenerate triangular
 dynamical $r$-matrix valued in $\wedge^2 \frakg_1$.
\end{enumerate}
\end{cor}
\pf By Proposition \ref{pro:split}, $T\frakh^* \times \{0\}$ is a 
Lie subalgebroid of $B$. Hence for any $X\in {\frak X} (\frakh^* ), 
(X, 0)\in \gm (B)$. Let $\phi_{t}$ be the (local) flow on $\frakh^* $
generated by $X$. The bisection 
$\exp{t(X, 0)}$  on the groupoid 
$\gm =\frakh^* \times \frakh^* \times G$ generated by the section 
$(X, 0)\in \gm (A)  $ 
is $\{(\lambda , \phi_{t}(\lambda ), 1 )|\lambda \in \frakh^* \}$.
Hence its   induced isomorphism between
$\gm_{\lambda}$ and $\gm_{\phi_{t}(\lambda )}$ 
is the  identity map, when  both   of them
are naturally  identified with $G$.  Here $\gm_{\lambda}$ and $\gm_{\phi_{t}(\lambda )}$
denote the isotropic groups of $\gm$ at the points
$\lambda$ and $\phi_{t}(\lambda )$, respectively. Therefore,
$Ad_{\exp{t(X, 0)} }$ is  an identity map between their
corresponding isotropic Lie algebras.
On the other hand, since $(X, 0)\in \gm (B)$,
hence $Ad_{\exp{t(X, 0)} }$, when being restricted
to $B$,  is exactly the map which establishes
the  isomorphism between $\frakg_{\lambda}$ 
and  $\frakg_{\phi_{t}(\lambda )}$.  Hence, $\frakg_{\lambda}$ and
 $\frakg_{\phi_{t}(\lambda )}$ are equal as Lie subalgebras of
$\frakg$.

For the second part, since $r$ is splittable, we have
$ r(\lambda )^{\#} \frakg^* \subset \frakg_1$ according to Proposition \ref{pro:split}.
Hence  $ \forall \lambda \in \frakh^* , \ r(\lambda )\in \wedge^{2}(r (\lambda )^{\#}
 \frakg^* )\subset \wedge^{2} \frakg_1$.
By dimension counting, one  easily sees that $r$ is non-degenerate when
being considered as a dynamical $r$-matrix  valued in $\wedge^2  \frakg_1$. \qed

Let $g: \frakh^* \lon G^{H}$ be a smooth map, where $G^{H}$ 
denotes the centralizer of $H$ in $G$ with its  Lie algebra being
denoted by $\frakg^{H}$.  Then $g$ can be naturally considered as a bisection
of the groupoid $\gm =\frakh^* \times \frakh^* \times G$, and hence
 we can talk about the induced automorphism $Ad_{g}$ of the
corresponding Lie algebroid. In particular, we have
a  Gerstenhaber algebra automorphism $Ad_{g}$ on 
$ \oplus \gm (\wedge^{*}A) $  \cite{Xu3}.

Given a smooth function $r:\frakh^* \lon \wedge^{2}\frakg$,
let $\Lambda_{r} =
\sum_{i} h_{i} \wedge \parr{}
 +r(\lambda )\in \gm ( \wedge^2  A)$ as in Proposition
\ref{pro:triangular}. Then
\be 
Ad_{g}\Lambda_{r} &=& Ad_{g}( \sum_{i} h_{i}\wedge \parr{}+r)\\
&=&\sum_{i}  Ad_{g}h_{i}\wedge ( \parr{}- \parr{g} g^{-1})+Ad_{g}r\\
&=&\sum_{i}  h_{i}\wedge ( \parr{}- \parr{g} g^{-1})+Ad_{g}r\\
&=&\sum_{i}  h_{i}\wedge\parr{} + (Ad_{g}r-\sum_{i} 
 h_{i}\wedge  \parr{g} g^{-1}).
\ee

Here in the second from the last equality, we used
$Ad_{g}h_{i}= h_{i}$ since $g\in G^H$.
Let 
\begin{equation}
r_{g}=Ad_{g}r-\sum_{i} h_{i}\wedge   \parr{g} g^{-1}.  
\end{equation}

Combining with Proposition \ref{pro:triangular},
 we thus have proved the following:

\begin{pro}
Assume  that  $g: \frakh^* \lon G^{H}$ is a smooth map. Then
\begin{enumerate}
\item $\Lambda_{r_g}=Ad_{g}\Lambda_{r} $;
\item  $r$ is a triangular  dynamical $r$-matrix iff $r_{g}$ is a triangular 
dynamical $r$-matrix.
\item \mbox{rank}$r_g$=\mbox{rank}$r$; in particular, if
$r$ is non-degenerate, so is $r_g$.
\end{enumerate}
\end{pro}

This proposition naturally  leads us to   the notion of gauge
 transformations on dynamical $r$-matrices, which was
 first introduced by Etingof and Varchenko \cite{EV1}.
Recall that triangular  dynamical $r$-matrices $r_1$ and $r_2$ are said to
be {\em gauge equivalent} if there exists a smooth
function $g:\frakh^* \lon G^H$ such that $r_2 =(r_{1})_{g}$. \\\\
{\bf Remark} Although non-degenerate  triangular 
   dynamical $r$-matrices are preserved by gauge transformations,
splittable dynamical $r$-matrices in general  are not. For example, the trivial
triangular dynamical  $r$-matrix $r=0$ is always splittable.
However $r_{g}= -\sum  h_{i}\wedge  \parr{g} g^{-1}$ is 
never  splittable unless $G^H =H$.\\\\

By $\calm (\frakg, \frakh )$, we denote the quotient space of the space
of all triangular  dynamical  $r$-matrices $r: \frakh^* \lon \wedge^2 \frakg$
 by  gauge transformations,   which is called
the {\em  moduli  space} of   triangular  dynamical  $r$-matrices.

Next we will introduce the dynamical  $r$-matrix  cohomology $H^*_{r} (\frakg ,
\frakh)$, whose second  cohomology group 
describes the tangent space of the moduli space
 $\calm (\frakg, \frakh )$. As we will see in Section 6,
the second  cohomology group $H^{2}_{r}(\frakg, \frakh )$
 is connected with the classification
of quantizations of  $r$ when it  
is  non-degenerate.

Consider $C^{k} =C^{\infty}(\frakh^* , (\wedge^k \frakg )^{H})$
(or equivalently denoted as $C^{\infty}(\frakh^* , (\wedge^k \frakg )^{\frakh})$),
 and define a differential  $\delta_{r} :C^{k}\lon C^{k+1}$ by
\begin{equation}
\label{eq:deltar}
\delta_{r} \tau =\sum_{i}h_{i}\wedge \parr{\tau }+[r , \tau ], \ \ \ \forall \tau 
\in C^{k}. 
\end{equation}

\begin{pro}
\label{pro:delta}
$\delta_{r}: 
C^{k}\lon C^{k+1}$ is well-defined and $\delta^2_{r} =0$.
\end{pro}
\pf It is clear that $\delta_{r} \tau$ is 
in $C^{\infty}(\frakh^* , (\wedge^{k+1} \frakg )^{H})$  provided that
$\tau \in C^{\infty}(\frakh^* , (\wedge^k \frakg )^{H})$.
 For any $\tau \in C^{k}=C^{\infty}(\frakh^* , (\wedge^k \frakg )^{H})$,
$\tau$ can be naturally considered as a section of
 $\wedge^{k}A$,  and
 \be
&&[\Lambda , \tau ]\\
&=&[\sum_{i} h_{i}\wedge \parr{} +r , \ \tau ]\\
&=&\sum_{i}h_{i}\wedge \parr{\tau }+[r , \tau ]\\
&=&\delta_{r} \tau.
\ee
Since $[\Lambda , \Lambda ]=0$,  it thus follows that
$\delta_{r}^2 =0$. \qed

Hence the  cochain complex $\delta_{r} : C^{k}\lon C^{k+1}$ defines a
cohomology, called {\em the dynamical $r$-matrix cohomology},
and denoted by $H_{r}^{*}(\frakg, \frakh )$. Two remarks are
in order.\\\\
{\bf Remark} \  (1). 
The cochain complex $\delta_{r} : C^{k}\lon C^{k+1}$ is in fact
a subcomplex of the Lie algebroid cohomology cochain complex
$d_{*} : \gm (\wedge^{k}A)\lon \gm (\wedge^{k+1}A), \ d_{*}X=[\Lambda , X]$.
  Therefore it is easy to see that such a cochain complex
is always  defined for an arbitrary
dynamical $r$-matrix, which is not necessary triangular.

(2). When $r$ is  triangular, $H_{r}^{*}(\frakg, \frakh )$ can be naturally
identified with  a ``special" $G\times H$-invariant Poisson cohomology
of the Poisson manifold $(M, \pi )$, i.e., the 
cohomology obtained by restricting  the Poisson
cochain complex to $G\times H$-invariant multi-vector
fields  tangent to the fibers of
the fibration: $\frakh^* \times G\lon  \frakh^*$.

\begin{pro}
\label{pro:delta1}
If  $g: \frakh^* \lon G^{H}$ is a smooth map, then
\begin{enumerate}
\item $\delta_{r_g}\smalcirc Ad_{g}=Ad_{g}\smalcirc \delta_{r}$;
\item $Ad_{g}: (C^{*} , \delta_{r})\lon (C^{*} , \delta_{r_g})$
induces an isomorphism $H^*_{r}(\frakg, \frakh )\cong H^*_{r_g}
(\frakg, \frakh )$.
\end{enumerate}
\end{pro}
\pf For any $\tau\in C^{\infty}(\frakh^* , (\wedge^{k}\frakg)^{H})$,
\be
  (Ad_{g}\smalcirc \delta_{r})\tau  &=& Ad_{g}[\Lambda , \tau ]\\
&=& [Ad_{g} \Lambda ,  Ad_{g} \tau ]\\
&=& [\Lambda_{r_{g}},  Ad_{g} \tau ] \\
&=& (\delta_{r_{g}}\smalcirc  Ad_{g})\tau. 
\ee
The conclusion thus follows immediately.
\qed

As a consequence, we conclude that
 $H_{r}^{*}(\frakg, \frakh )$ only depends
on the gauge equivalence class of the  dynamical $r$-matrix.
 For this reason,  we also  denote this group by
 $H_{[r]}^{*}(\frakg, \frakh )$.

\begin{pro}
\label{pro:tangent}
For any triangular dynamical $r$-matrix $r: \frakh^* \lon \wedge^{2}\frakg$,
$T_{[r]}\calm  (\frakg, \frakh )\cong  H^{2}_{[r]}(\frakg, \frakh )$.
\end{pro}
\pf In Equation (\ref{eq:cdybe}), 
 replace $r$ by $r+t\tau$ and take the derivative
at $t=0$, one obtains the linearization equation: 
$\sum_{i}h_{i}\wedge \parr{\tau }+[r , \tau ]=0$, i.e.,
 $\delta_{r} \tau =0$. It is clear that $\tau $ is
of  zero weight since $r+t\tau $ is of  zero weight.

To compute  the tangent space to the gauge orbit at $r$, one
needs to  compute $\frac{d}{dt}|_{t=0}(r_{\exp{tf}} )$, for
$f\in C^{\infty}(\frakh^* , \frakg^{H})$.
Now $r_{\exp{tf}}=Ad_{\exp{tf}}r-\sum_{i}h_{i}\wedge \parr{\exp{tf}}(\exp{tf})^{-1}$. It is thus simple to see that $\frac{d}{dt}|_{t=0}(r_{\exp{tf}})=
[f, r]-\sum_{i}h_{i}\wedge \parr{f}=-\delta_{r} f$.
The conclusion thus follows immediately. \qed 

Given a Lie algebra $\frakg$, one may also
consider  classical triangular  dynamical $r$-matrices 
$r_{\hbar}: \frakh^* \lon
 (\wedge^2  \frakg  ) \flb\hbar \frb$ valued in
$\frakg \flb\hbar \frb$  such that
$r_{\hbar} (\lambda )=r (\lambda )+\hbar r_{1} (\lambda )+\cdots$.
 The gauge transformation
can be defined formally in an obvious way. Thus
one  can form the moduli  space $\calm  (\frakg \flb\hbar \frb, \frakh )$.
Assume that $r: \frakh^* \lon  \wedge^{2} \frakg $
is a triangular classical dynamical $r$-matrix.
From Proposition \ref{pro:tangent}, it follows that
$T_{[r]}\calm  (\frakg \flb\hbar \frb, \frakh )\cong  H^{2}_{[r]}(\frakg, \frakh )
\flb\hbar \frb$.
By {\em a formal neighbourhood  of $r$} in $\calm  (\frakg \flb\hbar \frb, \frakh )$,
denoted by $\calm_{r}  (\frakg \flb\hbar \frb, \frakh )$,
we mean  the subset in $\calm  (\frakg \flb\hbar \frb, \frakh )$
consisting of the classes of those elements $r+O(\hbar )$.
Then $H^{2}_{[r]}(\frakg, \frakh )
\flb\hbar \frb$ can be considered  as a   linearization of
$\calm_{r}  (\frakg \flb\hbar \frb, \frakh )$. In general,
these two spaces are different. However, when $r$ is 
non-degenerate,
they expect to be isomorphic,  which should follow from Moser lemma.

 In fact, as we will see 
in the  next  theorem,  when $r$ is non-degenerate, 
$H_{[r]}^{*}(\frakg, \frakh )$ is  isomorphic to    the
relative Lie algebra cohomology.

\begin{thm}
If $r: \frakh^* \lon \wedge^2 \frakg $
 is a non-degenerate dynamical $r$-matrix, 
then  $H_{[r]}^{*}(\frakg, \frakh )$ is isomorphic 
to $H^* (\frakg, \frakh )$, the relative Lie algebra
cohomology of the pair  $(\frakg , \frakh)$.
\end{thm}
\pf Since $r$ is non-degenerate, $(M, \pi)$ is a 
symplectic manifold. As it is well known,
 $\pi^{\#}: \Omega^{*}(M) \lon {\frak X}^{*} (M)$ induces an isomorphism
between the
 de Rham cohomology cochain complex and 
the Poisson cohomology cochain complex. 
Now a $k$-mutivector field $P \in {\frak X}^{k}(M) $ is in
$C^{k}$ iff (i) $P$ is left $G$-invariant and right $H$-invariant; 
and (ii) $d\lambda^{i}\per P=0, \ \forall i=1, \cdots ,l$.
This, however, is equivalent to that (i) $(\pi^{\#})^{-1}P$ is
both left $G$-invariant and right $H$-invariant; and  
(ii) $\vh{i}\per (\pi^{\#})^{-1}P=0$,  because 
$\pi^{\#}(d\lambda^{i})=\vh{i}$, $\forall i=1, \cdots ,l$,
and  $\pi$ is $G\times H$-invariant.
 Note that a
 $k$-form $\omega \in \Omega^{k}(M)$  is $H$-invariant 
and  satisfies $\vh{i}\per \omega =0, \ \forall i=1, \cdots , l$,
iff $\omega$ is the pull back of a $k$-form on the quotient space
$M/H$, i.e, $\omega=p^* \omega '$, where
$p: M\lon M/H$  is the projection and
 $\omega '\in \Omega^{k}(M/H)$.  Moreover, $\omega$ is left $G$-invariant iff $\omega '$ is
left $G$-invariant since the left $G$-action on $M$ commutes 
with the right $H$-action. In summary, we have
proved that the space $ (\pi^{\#})^{-1}(C^{k})$ 
can be naturally identified with the space of left $G$-invariant
$k$-forms on $M/H\cong \frakh^* \times G/H$. Under
such an identification,  the
differential  $\delta_{r}$ goes to the de-Rham differential.
Hence $H^{k}_{[r]}(\frakg, \frakh )$ is isomorphic to
the invariant de-Rham cohomology $H^{k}(\frakh^* \times G/H )^{G}$.
Since $G$ does not act on the first factor $\frakh^* $,
the latter is isomorphic to $H^{k}(G/H )^{G}$, which is in turn
isomorphic to the relative Lie algebra cohomology 
$H^{k}(\frakg, \frakh )$ \cite{BW}. \qed

\section{Quantization and  star products}

In this section, we investigate 
the relation between quantizations of a triangular dynamical
$r$-matrix and  star products on its associated
Poisson manifold $(M, \pi )$. The main theme is
to show that quantizing $r$ is equivalent to finding
a  certain  special type of star products on $M$.
Let us first introduce the precise definition of a quantization.

\begin{defi}
Let $r: \frakh^* \lon \wedge^{2}\frakg$ be a 
triangular dynamical $r$-matrix.  A quantization of $r$ is an element
$F(\lambda ) =1+\hbar F_{1}(\lambda ) +\oh \in
C^{\infty}(\frakh^{*} , U\frakg \otimes U\frakg ) \flb\hbar \frb$ satisfying
\begin{enumerate}
\item the zero weight condition:
 $[1\ot h + h\ot 1, \ F(\lambda )]=0, \ \forall h\in \frakh$;
\item the   shifted cocycle condition:
\begin{equation}
\label{eq:shifted0}
(\Delta  \ot  id )F(\lambda )  F^{12} (\lambda - \half \hbar h^{(3)})
 =  (id \ot  \Delta )  F (\lambda )F^{23}(\lambda +\half \hbar h^{(1)}  );
\end{equation}
\item the  normal condition:
\begin{equation}
\label{eq:co0}
 (\epsilon \ot id) F(\lambda )  =  1; \ \
(id \ot \epsilon) F (\lambda ) =  1; \ \ \mbox{and}
\end{equation}
\item the quantization condition:
$F_{1}^{12} (\lambda )-F_{1}^{21} (\lambda )=r(\lambda)$,
\end{enumerate}
where $\Delta: U\frakg \lon U\frakg \ot U\frakg $ is
the standard comultiplication,  $\epsilon:
 \ug \lon \complex$  is the counit map,
and 
$F^{12} (\lambda -\half \hbar h^{(3)}), \ F^{23} (\lambda +\half \hbar h^{(1)})$
are  $ \ug\ot \ug \ot \ug $-valued
functions on $\frakh^*$ defined by
\begin{eqnarray}
F^{12} (\lambda -\half \hbar h^{(3)}) &=&F(\lambda  )\ot 1-
\frac{\hbar}{2} \sum_{i} \parr{F}\ot h_{i}
+\frac{1}{2!}(-\frac{\hbar}{2})^{2}\sum_{i_{1}i_{2}}
\frac{\partial^{2}F}{\partial \lambda^{i_{1}}\partial \lambda^{i_{2}}}\ot h_{i_
{
1}}h_{i_{2}} \nonumber\\
&& \ \ \ \ +\cdots
+\frac{1}{k!} (-\frac{\hbar}{2})^{k}
\sum \frac{\partial^{k} F}{\partial \lambda^{i_{1}}
\cdots \partial \lambda^{i_{k}}} \ot h_{i_{1}}\cdots h_{i_{k}} +\cdots,
\label{eq:Fshifted}
\end{eqnarray}
 and  similarly for $F^{23} (\lambda +\half \hbar h^{(1)})$.
\end{defi}

The relation between this definition  of quantizations
and  the well known quantum dynamical Yang-Baxter equation (QDYBE) is 
explained  by  the following proposition, which can be proved
by a straightforward verification.

\begin{pro}
If $F(\lambda )$ is a quantization of a triangular
dynamical $r$-matrix $r(\lambda): \frakh^* \lon \wedge^{2}\frakg$,
then $R(\lambda )= F^{21}(\lambda )^{-1}F^{12}(\lambda )$ 
can be written as   $R(\lambda )= 1+\hbar r(\lambda )+\oh$ and
satisfies the quantum dynamical Yang-Baxter equation (QDYBE):
\begin{equation}
\label{eq:qdybe}
R^{12}(\lambda-\half \hbar h^{(3)})R^{13}(\lambda
+\half \hbar h^{(2)} )R^{23}(\lambda-\half\hbar h^{(1)})
=R^{23}(\lambda +\half\hbar h^{(1)} )R^{13}(\lambda -\half
\hbar h^{(2)})R^{12}(\lambda +\half \hbar h^{(3)}).
\end{equation}
\end{pro}
{\bf Remark} This is a symmetrized version of QDYBE, which is 
known \cite{ESS} to be equivalent to the non-symmetrized QDYBE:
$$R^{12}(\lambda+\hbar h^{(3)})R^{13}(\lambda )R^{23}(\lambda+\hbar h^{(1)})
=R^{23}(\lambda )R^{13}(\lambda +\hbar h^{(2)})R^{12}(\lambda ).$$
The reason for us to  choose the symmetrized QDYBE 
in this  paper is because
 it is related to the Weyl quantization, while the
 non-symmetrized QDYBE
is related to  the normal ordering  quantization, as indicated in
\cite{Xu:gpoid}. Since we will use Fedosov method
later on, the  Weyl quantization is obviously of some advantage. \\\\

To proceed, we need some preparation on notations.
Let $\A=\cald \ot \ugh$, where $\cald$ is  the algebra
of smooth differential operators on $\frakh^{*}$. Then $\cald \ot \ug$
can be naturally identified with the algebra of left $G$-invariant
differential operators  on $M$. 
 Hence $\A$ becomes a 
Hopf algebroid \cite{Xu:gpoid}  with base algebra
$R= C^{\infty}(\frakh^* ) \flb\hbar \frb$. The comultiplication
$$\Delta : \ \ \A\lon \A\otr \A\cong 
\cald\ot_{C^{\infty}(\frakh^* )} \cald
\ot \ug \ot \ugh$$
is a natural extension of the comultiplications on $\cald$ and on  $\ug$: 
$$\Delta (D\ot u)=\Delta D\ot \Delta u, \ \ \ \ \forall D\in \cald,
\ \ \mbox{ and } u\in \ug, $$
where $\Delta D$ is  the bidifferential operator on $\frakh^*$ given by
$(\Delta D)(f, g)=D(fg), \ \forall f , \ g\in C^{\infty}(\frakh^* )$
 and
$\Delta u \in \ug\ot \ug$ is the usual comultiplication on $\ug$.
 Let us   fix a basis in  $\frakh$, say $\{h_{1},  \cdots ,h_{l}\}$,
and let $\{\xi_{1}, \cdots , \xi_{l}\}$ be its dual basis,
which in turn
  defines a coordinate system $(\lambda^{1}, \cdots ,\lambda^{l})$
on $\frakh^{*}$.
 
Set
\begin{equation}
\label{eq:theta}
\theta = \half \sum_{i=1}^{l}  (h_{i}\ot \frac{\partial}{\partial \lambda^{i}}-
\frac{\partial}{\partial \lambda^{i}}\ot h_{i})
\in \A\ot \A , \ \mbox{ and } \Theta =\exp{\hbar \theta} \in \A\ot \A.
\end{equation}
Note that $\theta $, and hence $\Theta$, is independent   of
the choice of a  basis in  $\frakh $.

For each $D\in \cald \ot \ug$,
we denote by $\Vec{D}$ its corresponding  left $G$-invariant
differential operator on $M=\frakh^* \times G$.  We also use a similar
notation to denote multi-differential operators on $M$ as well.
Now let $r(\lambda): \frakh^* \lon \wedge^{2}\frakg$ be a
triangular dynamical $r$-matrix, and $M=\frakh^* \times G$
its associated  (regular) Poisson manifold with
Poisson tensor $\pi =\sum_{i} \vh{i}\wedge \parr{} +\Vec{r(\lambda )}$. 
It is simple to see that the Poisson brackets on $C^{\infty}(M)$
can be described as follows:
\begin{enumerate}
\item for any $f, g\in
C^{\infty}(\frakh^* )$,
$\{f, g \}=0$;
\item for any $f \in C^{\infty}(\frakh^* )$ and $g\in C^{\infty}(G)$,
$\{f, g \}=-\sum_{i}\frac{\partial f}{\partial \lambda^{i}}\vh{i}g$;
\item for any $f, g \in C^{\infty}(G)$, $\{f, g \}=
\Vec{r(\lambda )}(f, g)$.
\end{enumerate}
This Poisson bracket relation naturally  motivates  the 
following theorem, which is  indeed
the main theorem of this section.

\begin{thm}
\label{thm:star}
Let $(M, \pi )$ be the Poisson
manifold associated to a  triangular dynamical $r$-matrix as in Proposition 
\ref{pro:Poisson}.
Assume that $*_{\hbar}$ is a $G\times H$-invariant
star product on $(M, \pi )$ satisfying the properties:
\begin{enumerate}
\item for any $f, g\in
C^{\infty}(\frakh^* )$,
$$f(\lambda )*_{\hbar}g (\lambda )=f(\lambda )
 g(\lambda );$$
\item  for any $f (\lambda )\in C^{\infty}(\frakh^* )$
and $  g  (x)\in  C^{\infty}(G)$,
\be
f(\lambda )*_{\hbar}g (x)
&=& \vT (f, g)=
\sum_{k=0}^{\infty} (-\frac{\hbar}{2})^{k}
\frac{1}{k!} 
\frac{\partial^{k} f}{\partial \lambda^{i_{1}}\cdots \partial \lambda^{i_{k}}}
\vh{i_{1}}\cdots \vh{i_{k}}g, \\
g (x)*_{\hbar}f(\lambda )
&=&\vT (g, f)=\sum_{k=0}^{\infty} (\frac{\hbar}{2})^{k}
\frac{1}{k!} \vh{i_{1}}\cdots \vh{i_{k}}g
\frac{\partial^{k} f}{\partial \lambda^{i_{1}}\cdots \partial \lambda^{i_{k}}};
\ee
\item there is a smooth map $F: \frakh^* \lon U\frakg \ot U\frakg \flb\hbar \frb$
such that for any $f(x), \ g(x)\in  C^{\infty}(G)$,
\begin{equation}
f*_{\hbar}g= \Vec{F(\lambda )}(f, g).
\end{equation}
\end{enumerate}
Then $ F(\lambda )$ is a quantization of the dynamical $r$-matrix $r(\lambda)$.
Conversely, any quantization of $r(\lambda)$
corresponds to a  $G\times H$-invariant star  product on $M$
satisfying   the  properties (i)-(iii).
\end{thm}

A   $G\times H$-invariant star product 
on $M$ with properties (i)-(iii) is called  a {\em compatible}
star product. In other words,  Theorem \ref{thm:star}
can   be stated that a  quantization of $r(\lambda)$ is
equivalent to a compatible star-product on $M$.


To prove Theorem \ref{thm:star}, we need several lemmas.

\begin{lem}
\label{lem:theta}
$\Theta $ satisfies  the equation:
\begin{equation}
\label{eq:Theta}
[(\Delta \ot id )\Theta ] \Theta^{12} =[(id \ot  \Delta) \Theta ]
 \Theta^{23}  \ \ \mbox{ in } \A\ot \A\ot \A.
\end{equation}
\end{lem}
\pf Note that both sides of Equation (\ref{eq:Theta})  normally
are  elements in $\A\otr \A\otr \A$. In our situation, however,
they  indeed can be considered as elements in $ \A\ot \A\ot \A$.

Now
\be
&&[(\Delta \ot  id )\Theta ]  \Theta^{12}\\
&=&[(\Delta \ot  id ) \exp{\hbar \theta} ]  \exp{\hbar \theta^{12}}\\
&=&\exp\hbar [(\Delta \ot  id )\theta + \theta^{12}]\\
&=&\exp\half \hbar \sum_{i=1}^{k}(
h_{i}\ot 1\ot \parr{} +1\ot h_{i}\ot \parr{}
+h_{i}\ot \parr{} \ot 1 - \parr{}\ot 1\ot h_{i}-1\ot \parr{}\ot h_{i}
-\parr{}\ot h_{i}\ot 1) .
\ee
Here in  the second equality we
  used the fact that $(\Delta \ot  id ) \theta$ and
$\theta^{12}$ commute in  $\A\ot \A \ot \A$.
 
A similar  computation leads to the same expression for
$[(id \ot  \Delta) \Theta ] \Theta^{23}$.
This proves Equation (\ref{eq:Theta}). \qed

\begin{lem}
\label{lem:f23}
$\forall D_{1}, D_{2}, D_3 \in \A$,
and $\forall  f_{1} (\lambda ) \in   C^{\infty}(\frakh^* )$,
$f_{2}(x)\in  C^{\infty}(G )$, and $g(\lambda , x)\in
C^{\infty}(\frakh^* \times G)$,
$$\Vec{[(\Delta \ot id )F(\lambda )(D_{1}\ot D_{2} \ot D_3 )]}
(f_{1} (\lambda ), \ f_{2}(x), \ g(\lambda ,x ))
=\Vec{[F^{23}(\lambda )(D_{1}\ot D_{2} \ot D_3 )]}
(f_{1} (\lambda ), \ f_{2}(x), \ g(\lambda , x )) .$$
\end{lem}
\pf Write $F (\lambda )=\sum a_{\alpha \beta} (\lambda ) u_\alpha
\ot u_{\beta}$, with $u_\alpha ,  u_{\beta}\in U\frakg $ and
$a_{\alpha \beta} (\lambda )\in C^{\infty}(\frakh^* )\flb\hbar \frb$.
Then $$ ((\Delta \ot id )F(\lambda ))(D_{1}\ot D_{2} \ot D_3 )
=\sum a_{\alpha \beta} (\lambda ) \Delta u_\alpha
(D_{1} \ot D_2 ) \ot u_{\beta} D_{3}.$$
Hence
\be
&&\Vec{[(\Delta \ot id )F(\lambda )(D_{1}\ot D_{2} \ot D_3 )]}
(f_{1} (\lambda ), \ f_{2}(x), \ g(\lambda ,x ))\\
&=&\sum a_{\alpha \beta} (\lambda ) \Vec{\Delta u_\alpha
(D_{1} \ot D_2 )}(f_{1} (\lambda ), \ f_{2}(x))
(\Vec{u_{\beta} D_{3}}g)(\lambda ,x )\\
&=&\sum a_{\alpha \beta} (\lambda ) \Vec{u_\alpha} 
[(\Vec{D_{1}}f_{1}) (\lambda )
(\Vec{D_{2}}f_{2} )(x ) ](\Vec{u_{\beta} D_{3}}g)(\lambda ,x )\\
&=&\sum a_{\alpha \beta} (\lambda ) (\Vec{D_{1}}f_{1}) (\lambda )
 ((\Vec{u_{\alpha} D_{2}})f_{2}) (x ) (\Vec{u_{\beta} D_{3}}g)(\lambda ,x )\\
&=&\Vec{ D_{1}\ot F(\lambda )(D_{2}\ot D_{3})}(f_{1} (\lambda ), \ f_{2}(x), \ g(\lambda ,x))\\
&=&\Vec{F^{23}(\lambda )(D_{1}\ot D_{2} \ot D_3 )}
(f_{1} (\lambda ), \ f_{2}(x), \ g(\lambda , x )).  \ee
\qed

\begin{cor}
\label{cor:exchange}
$\forall f_{1} (\lambda )\in C^{\infty}(\frakh^* ), \ f_{2} (x)\in
C^{\infty}(G)$ and $ g(\lambda , x )\in  C^{\infty}(\frakh^* \times G)$,
$$\vFT (f_{1} (\lambda )*_{\hbar} f_{2} (x), \ \  g(\lambda , x ) )
=\vT (f_{1} (\lambda ), \ \ \vFT (f_{2} (x), \ g(\lambda , x ) )).$$
\end{cor}
\pf
\be
&&\vFT  (f_{1} (\lambda )*_{\hbar} f_{2} (x), \ \  g(\lambda ,x ) )\\
&=&\vFT  (\vT (f_{1} (\lambda ), \ f_{2} (x)), \ g(\lambda , x ) )\\
&=&\Vec{(\Delta \ot id )(F(\lambda )\Theta ) \Theta^{12}}
(f_{1} (\lambda ), \ f_{2}(x), \ g(\lambda , x ))\\
&=&\Vec{ (\Delta \ot id )F(\lambda ) (\Delta \ot id )\Theta  \Theta^{12}}
(f_{1}(\lambda ), \ f_{2}(x), \ g(\lambda , x )) 
\mbox{ (by Lemma \ref{lem:theta})}\\
&=&\Vec{ (\Delta \ot id )F(\lambda ) (id \ot  \Delta) \Theta \Theta^{23} }
(f_{1}(\lambda ), \ f_{2}(x), \ g(\lambda , x )) 
\mbox{ (by Lemma \ref{lem:f23})}\\
&=&\Vec{ F^{23}(\lambda ) (id \ot  \Delta) \Theta \Theta^{23} }
(f_{1}(\lambda ), \ f_{2}(x), \ g(\lambda , x )).
\ee

Let us write $\Theta =\sum D_{\alpha}\ot D_{\beta }$. Then
$(id \ot  \Delta) \Theta =\sum  D_{\alpha}\ot \Delta D_{\beta }$,
and
\be
&&\Vec{ F^{23}(\lambda ) (id \ot  \Delta) \Theta \Theta^{23} }
(f_{1}(\lambda ), \ f_{2}(x), \ g(\lambda , x ))\\
&=& \sum  \Vec{ [D_{\alpha}\ot F(\lambda ) \Delta D_{\beta } \Theta  ]}
(f_{1}(\lambda ), \ f_{2}(x), \ g(\lambda , x ))\\
&=& \sum (\Vec{ D_{\alpha}}f_{1})(\lambda )
 \Vec{F(\lambda ) \Delta D_{\beta } \Theta } (f_{2}(x), \ g(\lambda , x )).
\ee

Using the expansion $\Theta =\sum_{k=0}^{\infty}
(\frac{\hbar}{2})^{k} \frac{1}{k!} (\sum_{i=1}^{l}  
(h_{i}\ot \frac{\partial}{\partial \lambda_{i}}-
\frac{\partial}{\partial \lambda_{i}}\ot h_{i}) )^{k}$,
one obtains that

\be
&&\vFT  (f_{1} (\lambda )*_{\hbar} f_{2} (x), \ \  g(\lambda , x ) )\\
&=&\sum_{k=0}^{\infty} (-\frac{\hbar}{2})^{k}
\frac{1}{k!}
\frac{\partial^{k} f_{1}(\lambda )}{\partial \lambda^{i_{1}}\cdots \partial \lambda^{i_{k}}}
\Vec{F(\lambda )\Delta (h_{i_{1}}\cdots h_{i_{k}}) \Theta} ( f_{2} (x), \ \  g(\lambda , x ) )\\
&=&\sum_{k=0}^{\infty} (-\frac{\hbar}{2})^{k}
\frac{1}{k!}
\frac{\partial^{k} f_{1}(\lambda )}{\partial \lambda^{i_{1}}\cdots 
\partial \lambda^{i_{k}}} 
\Vec{\Delta (h_{i_{1}}\cdots h_{i_{k}} )  F(\lambda ) \Theta }
 ( f_{2} (x), \ \  g(\lambda , x ) )\\
&=&\sum_{k=0}^{\infty} (-\frac{\hbar}{2})^{k}
\frac{1}{k!}
\frac{\partial^{k} f_{1}(\lambda )}{\partial \lambda^{i_{1}}\cdots 
\partial \lambda^{i_{k}}} \vh{i_{1}}\cdots \vh{i_{k}}
[\vFT (  f_{2} (x), \  g(\lambda , x ))]\\
&=&\vT (f_{1} (\lambda ), \ \ \vFT (f_{2} (x), \ g(\lambda , x ) )).
\ee
Here  the second equality follows 
 from the  fact that $F(\lambda )$ is of zero weight, i.e.,
 $F(\lambda ) (\Delta h ) 
=(\Delta h )F(\lambda )$, 
 $\forall h\in \frakh$.  This concludes the proof. \qed

\begin{pro}
\label{pro:various}
Under the same hypothesis as in Theorem \ref{thm:star}, we have\\
(1). for any $f(\lambda )\in C^{\infty}(\frakh^* )$ and 
$g(\lambda , x)\in C^{\infty} (\frakh^* \times G )$, 
\begin{eqnarray}
f(\lambda ) *_{\hbar} g(\lambda , x)&=&
\vT (f, g)=
\sum_{k=0}^{\infty} (-\frac{\hbar}{2})^{k}
\frac{1}{k!}
\frac{\partial^{k} f}{\partial \lambda^{i_{1}}\cdots \partial \lambda^{i_{k}}}
\vh{i_{1}}\cdots \vh{i_{k}}g, \label{eq:7} \\
g (\lambda , x)*_{\hbar}f(\lambda )&=& \vT (g, f)=\sum_{k=0}^{\infty} (\frac{\hbar}{2})^{k}
\frac{1}{k!} \vh{i_{1}}\cdots \vh{i_{k}}g
\frac{\partial^{k} f}{\partial \lambda^{i_{1}}\cdots \partial \lambda^{i_{k}}};
\label{eq:8}
\end{eqnarray}
(2). for any $f(\lambda , x)\in C^{\infty} (\frakh^* \times G )$ and
$g(x )\in C^{\infty}(G )$,
\begin{eqnarray}
f(\lambda , x) *_{\hbar} g(x )&=&
(\vFT ) (f, g)=\sum_{k=0}^{\infty}  (-\frac{\hbar}{2})^{k}
\frac{1}{k!} \vF
(\frac{\partial^{k} f}{\partial \lambda^{i_{1}}\cdots \partial \lambda^{i_{k}}},
\   \ \vh{i_{1}}\cdots \vh{i_{k}}g), \label{eq:9} \\
g (x)*_{\hbar}f(\lambda , x)&=&(\vFT ) (g, f)
=\sum_{k=0}^{\infty} (\frac{\hbar}{2})^{k}
\frac{1}{k!}  \vF  (\vh{i_{1}}\cdots \vh{i_{k}}g, \  \ 
\frac{\partial^{k} f}{\partial \lambda^{i_{1}}\cdots \partial \lambda^{i_{k}}}).
\label{eq:10}
\end{eqnarray}
\end{pro}
\pf We will prove Equation (\ref{eq:7}) first.
For that,  it suffices to
show this  for $g(\lambda , x)=g_{1}(\lambda )*_{\hbar}g_{2}(x)$,
$\forall g_{1}(\lambda )\in C^{\infty}(\frakh^* )$ and $g_{2}(x) \in
C^{\infty}(G )$, since,  at each point, the $ C^{\infty}$-jet space of  
$C^{\infty} (\frakh^* \times G )\flb\hbar \frb$
is spanned by the $ C^{\infty}$-jets of this type of
functions. Now
\be
f(\lambda ) *_{\hbar} g(\lambda , x )&=&
f(\lambda ) *_{\hbar} (g_{1}(\lambda )*_{\hbar}g_{2}(x) )\\
&=&(f(\lambda ) *_{\hbar} g_{1}(\lambda ))*_{\hbar}g_{2}(x)\\
&=&(f(\lambda )  g_{1}(\lambda ))*_{\hbar}g_{2}(x)\\
&=&\vT (f(\lambda )  g_{1}(\lambda ), \ g_{2}(x) )\\
&=&\vT (\vT (f(\lambda ), \ g_{1}(\lambda )), \ g_{2}(x) )\\
&=&\vdiT (f(\lambda ), g_{1}(\lambda ), g_{2}(x) ) 
 \mbox{ (by Lemma \ref{lem:theta})}\\
&=&\vidT (f(\lambda ), g_{1}(\lambda ), g_{2}(x) )\\
&=& \vT (f(\lambda ), \ \vT (g_{1}(\lambda ), g_{2}(x)))\\
&=& \vT (f(\lambda ), \ g_{1}(\lambda )*_{\hbar}g_{2}(x) )\\
&=&  \vT (f(\lambda ), \ g(\lambda , x )).
\ee
Equation (\ref{eq:8}) can be proved similarly.

To prove Equation (\ref{eq:9}), similarly we may assume that
$f(\lambda , x)=f_{1} (\lambda )*_{\hbar} f_{2} (x)$, for
$f_{1} (\lambda )\in C^{\infty}(\frakh^* )$ and
$f_{2}(x)\in  C^{\infty}(G )$.
Then
\be
f(\lambda , x) *_{\hbar} g(x )&=&
(f_{1} (\lambda )*_{\hbar} f_{2} (x))*_{\hbar} g(x )\\
&=&f_{1} (\lambda )*_{\hbar} (f_{2} (x)*_{\hbar} g(x ) )
\mbox{ (using Equation (\ref{eq:7}))}\\
&=&\vT (f_{1} (\lambda ) , \ \ f_{2} (x)*_{\hbar} g(x ))\\
&=&\vT (f_{1} (\lambda ) , \ \vF (f_{2} (x), \ g(x )))\\
&=&\vT (f_{1} (\lambda ) , \ \vFT (f_{2} (x), \ g(x ))) 
\mbox{ (by Corollary \ref{cor:exchange})}\\
&=&\vFT  (f_{1} (\lambda )*_{\hbar} f_{2} (x), \ \  g( x ) )\\
&=&\vFT (f(\lambda , x), \ g( x ) ).
\ee
Equation (\ref{eq:10}) can also be proved similarly.  \qed

We are now ready to prove the main theorem of the section.

\begin{thm}
\label{thm:full-star}
Under the same hypothesis as in Theorem \ref{thm:star},
$\vFT$ is the formal bidifferential operator
defining the star product $*_{\hbar}$, i.e., for
any $f(\lambda , x), \ g(\lambda , x)\in C^{\infty} (\frakh^* \times G )$,
$$f(\lambda , x) *_{\hbar} g(\lambda , x)=\vFT ( f, \ g).$$
\end{thm}
\pf We may assume
that  $f(\lambda , x)=f_{1} (\lambda )*_{\hbar} f_{2} (x)$, for
$f_{1} (\lambda )\in C^{\infty}(\frakh^* )$ and
$f_{2}(x)\in  C^{\infty}(G )$.
Then
\be
&&f(\lambda , x) *_{\hbar} g(\lambda , x)\\
&=&(f_{1} (\lambda )*_{\hbar} f_{2} (x)) *_{\hbar} g(\lambda , x)\\
&=&f_{1} (\lambda )*_{\hbar} (f_{2} (x)*_{\hbar} g(\lambda , x))
\ \mbox{ (by Proposition \ref{pro:various})}\\
&=&\vT (f_{1} (\lambda ), \ \ \vFT (f_{2} (x), \ g(\lambda , x ) ) )
\ \mbox{ (by Corollary \ref{cor:exchange})}\\
&=& \vFT (f_{1} (\lambda )*_{\hbar} f_{2} (x), \ \  g(\lambda , x ) )\\
&=& \vFT (f(\lambda , x), \ \ g(\lambda , x )).
\ee
This concludes the proof. \qed

Finally,  before proving  Theorem  \ref{thm:star}, we need
the following result, which  connects the
shifted cocycle condition with  the associativity
of a  star-product.

\begin{pro}
\label{pro:f123}
Under the same hypothesis as in Theorem \ref{thm:star},
 $ \forall f_{1}(x), f_{2}(x) , f_{3}(x)\in C^{\infty}(G)$,
\begin{enumerate}
\item $\Vec{(\Delta  \ot  id )F(\lambda )  F^{12} (\lambda - \half
 \hbar h^{(3)})}(f_{1}(x),\ f_{2}(x) ,\ f_{3}(x))
=(f_{1}(x)*_{\hbar}f_{2}(x)) *_{\hbar}f_{3}(x)$;
\item $\Vec{ (id \ot  \Delta )
  F (\lambda )F^{23}(\lambda +\half \hbar h^{(1)}  )}(f_{1}(x),\ f_{2}(x) ,\ f_{3}(x))
=f_{1}(x)*_{\hbar}(f_{2}(x) *_{\hbar} f_{3}(x))$. 
\end{enumerate}
\end{pro}
\pf  From Equation (\ref{eq:Fshifted}), it follows that
$$(\Delta  \ot  id )F(\lambda )  F^{12} (\lambda - \half
 \hbar h^{(3)})= \sum \frac{1}{k!} (-\frac{\hbar}{2})^{k}
[(\Delta  \ot  id )F(\lambda ) ] (\frac{\partial^{k} F}{\partial \lambda^{i_{1}}
\cdots \partial \lambda^{i_{k}}} \ot h_{i_{1}}\cdots h_{i_{k}} ). $$
Hence
\be
&&\Vec{(\Delta  \ot  id )F(\lambda )  F^{12} (\lambda - \half
 \hbar h^{(3)})}(f_{1}(x),\ f_{2}(x) ,\ f_{3}(x))\\
&=&\sum \frac{1}{k!} (-\frac{\hbar}{2})^{k} \Vec{F(\lambda )}
[\Vec{\frac{\partial^{k} F(\lambda )}{\partial \lambda^{i_{1}}
\cdots \partial \lambda^{i_{k}}}} (f_{1}(x),\ f_{2}(x)), \ \ 
  (\vh{i_{1}}\cdots \vh{i_{k}} f_{3})(x)]\\
&=& \sum \frac{1}{k!} (-\frac{\hbar}{2})^{k} \Vec{F(\lambda )}
[\frac{\partial^{k}(f_{1} *_{\hbar}f_{2})}{\partial \lambda^{i_{1}}
\cdots \partial \lambda^{i_{k}}}, \ \   (\vh{i_{1}}\cdots \vh
{i_{k}} f_{3})(x)] \ (\mbox{using Equation (\ref{eq:9})})\\
&=& ( f_{1}(x)*_{\hbar}f_{2}(x)) *_{\hbar} f_{3}(x).
\ee
The second identity can be proved similarly. \qed
{\bf Proof of Theorem \ref{thm:star} }  Since $*_{\hbar}$ is
 invariant under the right $H$-action, $\Vec{F(\lambda )}$ is
right $H$-invariant. This implies that $F(\lambda )$ is
$Ad_{H}$-invariant,  and therefore is   of zero weight.
The normal condition
follows from the fact  that $1$ is the unit of the star algebra,
i.e., $1*_{\hbar}f=f*_{\hbar}1=f$.
And the shifted cocycle condition follows from the associativity
of  the star product together with  Proposition \ref{pro:f123}.
Finally,   let us write  $F(\lambda )=1+\hbar F_{1} (\lambda )+\oh$.
Since $*_{\hbar}$ is a star product quantizing  $\pi$,  it follows
that $\Vec{(F_{1} (\lambda )-F_{1}^{21} (\lambda ))} (f, g)= \{f, \ g\}=
\Vec{r(\lambda)} (f, g), \  \forall f , g\in C^{\infty}(G)$.
Hence it follows that $F_{1} (\lambda )-F_{1}^{21} (\lambda )=r(\lambda)$.

Conversely, if $F(\lambda )$ is a quantization 
of $r(\lambda)$, according to Theorem 7.5 in \cite{Xu:gpoid},
$\vFT$ is indeed  an associator and therefore  defines a 
star product on $M=\frakh^* \times G$. It is simple to
see that this star product  is a quantization of $\pi$ and
satisfies Properties (i)-(iii) in Theorem \ref{thm:star}.  \qed
We end this section by the following\\\\
{\bf Remark} Bordemann et. al.  found an explicit
formula for a star-product on
 $\reals \times SU (2)$ \cite{BBEW} using a 
quantum analogue of Marsden-Weinstein reduction. It would be
interesting to see if this is a compatible star-product.

\section{Symplectic connections}

From now on, we will confine ourselves mostly to non-degenerate
triangular dynamical $r$-matrices. In this case, the corresponding
Poisson manifolds  are  in fact symplectic, and
therefore can be quantized by Fedosov method \cite{F1, F2}. As is
well known, Fedosov quantization relies on
the choice of a symplectic connection. Serving
as a preliminary,
this section is devoted to the discussion on symplectic connections.
We will start with some general   notations and constructions.

Let $\nabla$ be a torsion-free symplectic connection on a symplectic
manifold $(M,  \omega )$.
 Define the symplectic curvature \cite{F1} by
\begin{equation}
\label{eq:sym-curvature}
R( X, Y, Z, W)=\omega (X, R(Z, W)Y), \ \ \forall X, Y, Z, W\in {\frak X}  (M),
\end{equation}
where $R(Z, W)Y=\nabla_{Z}\nabla_{W}Y -\nabla_{W}\nabla_{Z}Y-\nabla_{[Z, W]}Y$
is the usual curvature tensor of $\nabla$.

\begin{pro}
\label{pro:curvature}
\begin{enumerate}
\item $R( X, Y, Z, W)$ is skew symmetric with respect to $Z$ and $W$, and
symmetric with respect to $X$ and $ Y$, i.e.,
\begin{equation}
R( X, Y, Z, W)=-R( X, Y, W, Z), \ \ \ \ R( X, Y, Z, W)=R(Y, X, Z, W).
\end{equation}
\item The following Bianchi's identity holds:
\begin{equation}
R( X, Y, Z, W)+R(X, Z, W, Y)+R(X, W, Y, Z)=0.
\end{equation}
\end{enumerate}
\end{pro}
\pf  It is clear by definition that $R( X, Y, Z, W)$ is skew symmetric with 
respect to $Z$ and $W$. Now since  $\nabla$ is a
symplectic connection, then
\be
&&\omega (X, \nabla_{Z}\nabla_{W}Y )\\
&=&Z(\omega (X, \nabla_{W}Y ))-\omega (\nabla_{Z}X, \nabla_{W}Y )\\
&=&Z(W\omega (X, Y)) -Z\omega (\nabla_{W}X, Y)-W\omega (\nabla_{Z}X, Y)
+\omega (\nabla_{W} \nabla_{Z} X, Y).
\ee

Similarly,
$$\omega (X, \nabla_{W}\nabla_{Z} Y)=W(Z\omega (X, Y)) -W\omega (\nabla_{Z}X, Y)-Z\omega (\nabla_{W}X, Y)
+\omega (\nabla_{Z} \nabla_{W} X, Y).$$
Hence
$$\omega (X, \nabla_{[Z, W]}Y)=[Z, W](\omega (X, Y))-\omega (\nabla_{[Z, W]}X, Y).$$
Thus
\be
R( X, Y, Z, W)&=&\omega (X,\  R(Z, W)Y)\\
&=&\omega (X,  \ \nabla_{Z}\nabla_{W}Y -\nabla_{W}\nabla_{Z}Y-\nabla_{[Z, W]}Y)\\
&=&-\omega (\nabla_{Z}\nabla_{W}X -\nabla_{W}\nabla_{Z}X-\nabla_{[Z, W]}X, \ Y)\\
&=&\omega (Y, R(Z, W)X)\\
&=&R(Y, X, Z, W).
\ee
This concludes the proof of (i). Finally, (ii) follows from
the usual Bianchi's identity for a torsion-free connection.   \qed

Symplectic connections always exist on any symplectic
manifold. In fact, there is a standard procedure to construct
 a torsion-free symplectic connection  from 
 an arbitrary  torsion-free linear connection \cite{L, F1}.
Since such a   construction is essential  to our discussion here, 
let us recall it briefly below.

Assume that $\nablaa $ is a   torsion-free linear  connection  on a symplectic
manifold $M$. Then any   linear  connection on $M$ can be written
as
\begin{equation}
{\nabla}_{X}Y=\nablaa_{X}Y +S(X, Y) , \ \ \ \forall X, Y\in {\frak X}  (M),
\end{equation}
where $S$ is a $(2, 1)$-tensor on $M$.
Clearly, ${\nabla}$ is torsion-free iff
$S$ is symmetric, i.e., $S(X, Y)=  S(Y, X)$,  $\forall X, Y\in
 {\frak X} (M)$. And
${\nabla}$ is symplectic iff ${\nabla}_{X}\omega =0$.
The latter is equivalent to
\begin{equation}
\omega (S(X, Y), Z)-\omega (S(X, Z), Y)=(\nablaa_{X}\omega )(Y, Z),
\ \ \forall X, Y, Z\in {\frak X}  (M).
\end{equation}

\begin{lem} 
\label{lem:sym-conn}
If $\nablaa $ is a   torsion-free linear  connection,  and
 $S$  is a  $(2, 1)$-tensor defined by
the equation:
\begin{equation}
\label{eq:s}
\omega (S(X, Y), Z)=\frac{1}{3}[(\nablaa_{X}\omega )(Y, Z)+
(\nablaa_{Y}\omega )(X, Z)] ,
\end{equation}
then ${\nabla}_{X}Y=\nablaa_{X}Y +S(X, Y)$ is a torsion-free
symplectic connection. Moreover, if $M$ is a symplectic $G$-space
and $\nablaa$ is a $G$-invariant connection, then 
  ${\nabla}$ is also $G$-invariant.
\end{lem}
\pf Clearly, $S(X, Y)$,   defined  in this way,  is symmetric
with respect to $X$ and $Y$.  Now
\be
&&\omega (S(X, Y), Z)-\omega (S(X, Z), Y)\\
&=&\frac{1}{3}[(\nablaa_{X}\omega )(Y, Z)+ (\nablaa_{Y}\omega )(X, Z)]
-\frac{1}{3}[(\nablaa_{X}\omega )(Z, Y)+(\nablaa_{Z}\omega )(X, Y)]\\
&=&\frac{1}{3}[(\nablaa_{X}\omega )(Y, Z)+(\nablaa_{Y}\omega )(X, Z)
+(\nablaa_{X}\omega )( Y, Z)+(\nablaa_{Z}\omega )(Y, X)]\\
&=&(\nablaa_{X}\omega ) (Y, Z),
\ee
where the last step follows from  the identity:
 
$$(\nablaa_{X}\omega )(Y, Z)+(  \nablaa_{Y}\omega )(Z, X)+
(\nablaa_{Z}\omega )(X, Y)=0.$$
 
This means that ${\nabla}$ is a torsion-free symplectic connection. 
The second statement is obvious according to Equation (\ref{eq:s}).\qed

Now  we  retain to the case that  $M=\frakh^* \times G$,
the symplectic manifold associated with a non-degenerate
triangular dynamical $r$-matrix $r$, which is our  main
subject  of interest in the present paper. The main result
is the following

\begin{thm}
\label{thm:connection}
Assume that $r: \frakh^* \lon \wedge^2 \frakg$ is a
non-degenerate triangular dynamical $r$-matrix.
Let $M=\frakh^* \times G$ be equipped with 
the symplectic structure as in Corollary \ref{cor:sym}. Then
$M$ admits a $G\times H$-invariant torsion-free
 symplectic connection  $\nabla$ satisfying
the property that $\nabla_{X}\vh{}=0, \ \forall  X\in {\frak X}  (M),
\ h\in \frakh$.
\end{thm}

We need a couple of lemmas first.
\begin{lem}
\label{lem:del0}
Assume that $\frakg$ admits a reductive decomposition
$\frakg =\frakh \oplus \muu$, i.e., $[\frakh , \muu ]\subset
 \muu$. Then, the following equations define
 a  biinvariant torsion-free linear connection $\nablaa$
 on $M$:

\begin{eqnarray}
\begin{array}{lll}
\nablaa_{X}\lam{i}= 0, & \nablaa_{X} \Vec{h}= 0,&  \nablaa_{X} \Vec{e}=0;  \\
\nablaa_{\vh{}}\lam{i}=0,& \nablaa_{\vh{}}\vh{1}=0,& \nablaa_{\vh{}} \ve{} =
\Vec{[h, e]}; \\
\nablaa_{\ve{}}\lam{i}=0, & \nablaa_{\ve{}}\vh{} =0,& \nablaa_{\ve{1}}
\ve{2}= \half \Vec{[e_{1}, e_{2}]},
\end{array}
\end{eqnarray}
where  $X \in {\frak X}  (\frakh^* )$,
$h, h_{1}\in \frakh $,  and $e, \ e_{1},\  e_{2}\in  \muu$.
\end{lem}
\pf This follows from a straightforward verification. \qed 

\begin{lem}
\label{lem:decom}
Given a Lie algebra $\frakg$, if there exists a
non-degenerate triangular dynamical $r$-matrix
  $r: \frakh^* \lon \wedge^2 \frakg$, then $\frakg$
 admits a  reductive decomposition $\frakg =\frakh \oplus \muu$ so that
$[\frakh  , \muu ]\subset \muu$.
\end{lem}
\pf Fixing  any $\lambda \in \frakh^*$, we  take $\muu =r(\lambda )^{\#} 
\frakh^{\perr}$.  Since $r(\lambda )$ is
non-degenerate, by definition, we have $\frakg =\frakh +\muu $.
On the other hand, it is clear that $\mbox{dim} \muu\leq
\mbox{dim}\frakh^{\perr}=\mbox{dim}\frakg -\mbox{dim}\frakh$.
Hence, $\mbox{dim}\frakh +\mbox{dim} \muu 
\leq \mbox{dim}\frakg$. Therefore $\frakg =\frakh +\muu $
 must be a direct sum.  For any  
$h\in \frakh$ and $\xi\in  \frakg^*$, since $r(\lambda )$
is of zero weight,  we have  $[h, r(\lambda )^{\#}\xi ]
=r(\lambda )^{\#} (ad_{h}^* \xi  )$. Since $ad_{h}^* \xi\in \frakh^{\perr}$
for any $\xi \in \frakg^*$, it follows that $\muu =r(\lambda )^{\#} 
\frakh^{\perr}$ is stable under the adjoint action of $\frakh$. \qed
{\bf Remark} Note that, in our proof above,
the decomposition $\frakg =\frakh \oplus \muu $
depends on the choice of a  particular  point $\lambda \in \frakh^*$.
 It is not clear if $\muu =r(\lambda )^{\#} 
\frakh^{\perr} $ is independent of  $\lambda$.   \\\\\\
{\bf Proof of Theorem \ref{thm:connection}} 
 According to Lemma  \ref{lem:decom},
we may find a  reductive decomposition $\frakg=\frakh \oplus \muu$ such that
$[\frakh  , \muu ]\subset \muu$. Let $\nablaa$ be
the  $G$-biinvariant  torsion-free  connection on $M$ as
in  Lemma \ref{lem:del0}.  According to  Lemma \ref{lem:sym-conn},
 one can construct
 a torsion-free symplectic connection  $\nabla $ on $M$.
Since the symplectic structure is $G\times H$-invariant,
the resulting symplectic connection 
$\nabla $ is $G\times H$-invariant. It remains to
show that $\nabla $ still satisfies the condition
that $\nabla_{X}\vh{}=0, \ \forall  X\in {\frak X}  (M)$ and 
$ h\in \frakh$.  The latter is equivalent to that $S(\vh{}, X)=0$.
To show this identity, first note that
$\forall X\in {\frak X}  (M)$,  $\nablaa_{\vh{}}X=L_{\vh{}}X$, since
$\nablaa$ is torsion-free and $\nablaa_{X}\vh{}=0$. Hence
$\nablaa_{\vh{}}\omega =L_{\vh{}}\omega $. However, 
$L_{\vh{}}\omega =0$ since $\omega $ is invariant under
 the right $H$-action.
Thus, we have $\nablaa_{\vh{}}\omega=0$. According to Equation (\ref{eq:s}),
$\forall Y\in {\frak X}  (M)$, $\omega (S(\vh{}, X), Y)=
\third [(\nablaa_{\vh{}}\omega )(X, Y)+(\nablaa_{X}\omega )(\vh{}, Y)]
=\third (\nablaa_{X}\omega )(\vh{}, Y)$. This implies that
$\omega^{b }  (S(\vh{}, X))=\third (\vh{} \per \nablaa_{X}\omega )
=\third   \nablaa_{X} (\vh{}  \per \omega )$ since
$\nablaa_{X} \vh{}  =0$. Finally, for any $i$, $\vh{i}  \per \omega
=d\lambda^{i}$ and from the table in Lemma \ref{lem:del0}, it
is easy to check that $\nablaa_{X} (d\lambda^{i})=0, \ \forall i=1, \cdots , l$.
 It thus follows that $S(\vh{i}, X)=0,  \ \forall i=1, \cdots , l$.
This concludes the proof. \qed

In the case that $r(\lambda )\in \wedge^{2}\muu$, the symplectic
connection can be described more explicitly.
\begin{pro}
\label{pro:special}
Suppose that $\frakg =\frakh\oplus \muu$ is a reductive decomposition,
  $\{h_{1}, \cdots , h_{l}\}$
is a basis of $\frakh$,  and $\{e_{1}, \cdots , e_{m}\}$ is a  basis
of $\muu$. Suppose that $ r(\lambda )=\sum_{ij} r^{ij}(\lambda )e_{i}\wedge
e_{j}$ is a non-degenerate triangular dynamical $r$-matrix. Then
the symplectic   connection on $M$  obtained from $\nabla^0$, using
the standard construction as in Lemma \ref{lem:del0},  has the following
form:
\begin{eqnarray}
\begin{array}{lll}
\nabla_{\lam{i}}\lam{j}=0, & \nabla_{\lam{i}} \vh{j}= 0,& \nabla_{\lam{i}}\ve{j}
=\sum_{k}d_{ij}^{k}(\lambda ) \ve{k}; \\
\nabla_{\vh{i}}\lam{j}=0,& \nabla_{\vh{i}}\vh{j}=0,& \nabla_{\vh{i}} \ve{j} =
\Vec{[h_{i}, e_{j}]}; \\
\nabla_{\ve{i}}\lam{j}=\sum_{k}d_{ij}^{k}(\lambda ) \ve{k} ,&
 \nabla_{\ve{i}}\vh{j} =0, &\nabla_{\ve{i}}\ve{j}= \half
\Vec{[e_{i}, e_{j}]} + \sum_{k} f_{ij}^{k}(\lambda ) \ve{k},
\end{array}
\end{eqnarray}
where $d_{ij}^{k}(\lambda )$ and $f_{ij}^{k}(\lambda )$ are 
smooth functions on $\frakh^*$.
\end{pro}
\pf  The proof is essentially a straightforward
computation.  We omit it here. \qed

\begin{cor}
Under the same hypothesis as in Proposition \ref{pro:special}, if
$\{\hs{1}, \cdots , \hs{l}, \es{1}, \cdots , \es{m}\}$
denotes the dual basis of $\{h_{1},  \cdots , h_{l},
e_{1}, \cdots , e_m \}$, then
\begin{eqnarray}
\begin{array}{lll}
\nabla_{\lam{i}}d\lambda^{j}=0, & \nabla_{\lam{i}} \vhs{j}= 0,& 
\nabla_{\lam{i}}\ves{j}
=\sum_{k}d_{ik}^{j}(\lambda ) \ves{k}; \\
\nabla_{\vh{i}}d\lambda^{j}=0,& \nabla_{\vh{i}}\vhs{j}=0,& 
\nabla_{\vh{i}} \ves{j} = \Vec{ad_{h_{i}}^{*}\es{j}} ;\\
\nabla_{\ve{i}} d\lambda^{j}=0 ,&
 \nabla_{\ve{i}}\vhs{j} =-\half \sum_{k} a_{ik}^{j}\ves{k},
 &\nabla_{\ve{i}}\ves{j}=  -\sum_{k}d_{ik}^{j}(\lambda )
 d\lambda^{k} -( \half a_{ik}^{j}+f_{ik}^{j}(\lambda ))\ves{k},
\end{array}
\end{eqnarray}
where the coadjoint action is defined
by $<ad^{*}_{u}\xi , v>=-<\xi , [u, v]>$, 
 $\forall u, v\in \frakg$ and $\xi \in \frakg^*$, and
   the constants $a_{ij}^{k}$ are defined by the
equation $[e_{i}, e_{j}]=\sum_{k}a_{ij}^{k}h_{k} \ (\mbox{mod } \muu)$.
\end{cor}

We end this section by generalizing Theorem \ref{thm:connection}
to the  splittable triangular dynamical $r$-matrix case.
According to Corollary \ref{cor:split-nondeg}, one
may reduce a splittable triangular dynamical $r$-matrix
to a non-degenerate one by considering   the
Lie subalgebra $\frakg_{1}\subset \frakg$. Thus immediately we  obtain
the following

\begin{cor}
\label{cor:conn-Poisson}
Assume that $r: \frakh^* \lon \wedge^2 \frakg$ is a
splittable triangular dynamical $r$-matrix.
Let $M=\frakh^* \times G$ be  its associated Poisson manifold
as in Proposition \ref{pro:Poisson}, which admits a 
(regular) symplectic foliation.
Then there exists  a $G\times H$-invariant torsion-free
 leafwise  Poisson connection  $\nabla$ satisfying
 $\nabla_{X}\vh{}=0$, for any  $h\in \frakh$
and any vector field   $X\in {\frak X} (M)$   tangent to the symplectic
foliation.
\end{cor}

However, when a triangular dynamical $r$-matrix
 $r$ is not splittable, such a Poisson
connection may not exist. We give a counterexample
below. \\\\\
{\bf Example 4.9.} Consider a two dimensional Lie algebra
$\frakg$ with basis $\{h, e\}$ satisfying  the
bracket relation $[h, e]=ah$, where $a$ is a fixed constant.
 Let $\frakh =\reals h$ and
$r(\lambda )=f(\lambda )h\wedge e$, where $f(\lambda )$ is a
smooth function. It is simple to see that $r(\lambda )$ is
a triangular dynamical $r$-matrix of rank zero,
and it  is not splittable unless $a=0$. Nevertheless,
$r(\lambda )$ defines a regular rank 2  Poisson structure on
the three dimensional space $M=\reals\times G$ with
the Poisson tensor $\pi =\vh{}\wedge
\frac{d}{d \lambda } +f(\lambda )
\vh{}\wedge \ve{}$, where $G$ is a 2-dimensional Lie group
integrating the Lie algebra $\frakg$. It is simple to see that
the symplectic foliation of $M$ is spanned by the vector fields 
$\vh{}$ and $\frac{d}{d \lambda }+f(\lambda )\ve{}$. 
Let us denote $X=\frac{d}{d \lambda }+f(\lambda )\ve{}$.
Then,   we have $[\vh{}, \ X]=af(\lambda ) \vh{}$.
Now suppose that $\nabla$ is a $G\times H$-invariant torsion-free
leafwise Poisson connection on $M$ satisfying the
condition that $\nabla_{\vh{}} \vh{}=0$ and $\nabla_X \vh{}=0$.
Since $\nabla$ is torsion-free, it  follows that
$\nabla_{\vh{}}X=[\vh{}, \ X]=af(\lambda ) \vh{}$.
Assume that $\nabla_{X}X=b(\lambda, x )\vh{}+c(\lambda, x )X$,
where $b(\lambda, x    )$ and $c(\lambda , x   )$ are smooth functions on $M$.
Then, $\nabla_{X}\pi =\nabla_{X} (\vh{}\wedge X)
=\vh{}\wedge \nabla_{X} X=c(\lambda , x   )\vh{}\wedge X$.
Since $\nabla$ is a Poisson connection, it follows that $c(\lambda, x )=0$.
Finally, we still need to check that $\nabla$ is
$G\times H$-invariant. It is clear that $\nabla$ is
$G$-invariant iff the function $b(\lambda, x )$ is
independent of $x\in G$ (which will be denoted by $b(\lambda )$).  For it to be
invariant under the right $H$-action, one needs the following
condition:
$$\nabla_{[\vh{}, X]}X+\nabla_{X}[\vh{}, X]=[\vh{}, \ \nabla_{X}X]=
[\vh{},\  b(\lambda )\vh{}]=0.$$
It thus follows that $\nabla_{(a f(\lambda )\vh{} )}X+\nabla_{X}(a f(\lambda )\vh{}
)=0$, which implies that 
$f^{2}(\lambda ) a^2 \vh{} +a(-\frac{d f}{d \lambda})\vh{}=0$.
 Therefore, we arrive at the following equation
(under the assumption that $a\neq 0$):
\begin{equation}
\frac{d f}{d \lambda} =a f^{2}(\lambda ).
\end{equation}
In conclusion, we have proved that such a connection does not exist unless
$f(\lambda )$ is a solution of the above equation. It would be
interesting to find out what is the geometric meaning
of this equation. \\\\\\
{\bf Remark} Our  quantization method
does not work for  this particular example. It  is  thus very
natural to ask whether this dynamical $r$-matrix is still
quantizable.  Etingof and Nikshych recently has  given an affimative
 answer to this question
 using the so called vertex-IRF transformation
method \cite{EN}.  Their method indeed works for a large class
of  dynamical $r$-matrices called ``completely degenerate", which
somehow is opposite to the non-degenerate ones   considered
in this paper.  It would be very interesting to see whether one
could combine these two methods together to
completely solve  the quantization problem
for arbitary triangular dynamical $r$-matrices.

\section{Compatible Fedosov star products}

  In this section, we consider Fedosov star products on
a symplectic Hamiltonian $H$-space $M$, where $H$ is  an Abelian group.
For the reader's convenience,
we will give a brief account of  the   general
construction of Fedosov star products in Appendix.
Readers may refer to that section for  various
notations  and formulas that are  used here.
What is eventually relevant 
to our situation is the case when $M$ is 
the symplectic manifold $\frakh^* \times G$ corresponding
to a nondegenerate dynamical $r$-matrix. However,
we believe that our general   presentation  would
 be of its own interest. 
We can now state the  main result of this section.


\begin{thm}
\label{thm:Fedosov}
Let $H$ be an Abelian group and  $M$ a symplectic  Hamiltonian $H$-space
with an equivariant momentum map  $J: M\lon \frakh^*$.
 Assume that $J$ is a submersion, and  
 there exists a $H$-invariant symplectic connection $\nabla$ such that
 $\vh{}$ is parallel for any $h\in \frakh$, i.e.,
$\nabla_{X}\vh{}=0$, $\forall X\in {\frak X} (M)$.
Let $*_{\hbar}$ be the corresponding  Fedosov star product on $M$ with
Weyl curvature $\Omega =\omega +\hbar \omega_{1} +\cdots +\hbar^{i}\omega_{i} 
+\cdots \in Z^{2}(M)\flb\hbar \frb$,
which satisfies the condition that $i_{\vh{}}\omega_{i}=0, \  \forall i\geq 1,\ \forall
h\in \frakh $. 
  Then for any $f(\lambda )\in C^{\infty}(\frakh^* )$
and $g(x)\in  C^{\infty}(M)$, we have
\be
(J^*f )*_{\hbar}g  (x)
&=& \sum_{k=0}^{\infty} (-\frac{\hbar}{2})^{k}
\frac{1}{k!}
J^*(\frac{\partial^{k} f}{\partial \lambda^{i_{1}}\cdots \partial \lambda^{i_{k}}} )
\vh{i_{1}}\cdots \vh{i_{k}}g; \\
g (x)*_{\hbar} (J^*f )
&=&\sum_{k=0}^{\infty} (\frac{\hbar}{2})^{k}
\frac{1}{k!} (\vh{i_{1}}\cdots \vh{i_{k}}g )
J^* (\frac{\partial^{k} f}{\partial \lambda^{i_{1}}\cdots \partial \lambda^{i_{k}}}).
\ee
Here $\vh{}$ denotes the  corresponding Hamiltonian vector field on $M$ generated
by $h\in \frakh$.
\end{thm}
{\bf Remark}  From  Theorem \ref{thm:Fedosov}, it follows  that
$J^* :  C^{\infty}(\frakh^* )\flb\hbar \frb\lon  C^{\infty}(M) \flb\hbar \frb$ is
an algebra homomorphism, where $C^{\infty}(\frakh^* )\flb\hbar \frb$ is
equipped with pointwise multiplication. In other words,
$J^*$ is a quantum momentum map \cite{Xu0}. It would be interesting to
see how to generalize this result to the case 
 when $H$ is not Abelian \cite{Xu:4}. \\\\\\

Applying Theorem \ref{thm:Fedosov} to the symplectic  manifold
$M=\frakh^* \times G$ associated  to a nondegenerate 
triangular dynamical $r$-matrix,
and using Theorem \ref{thm:connection}, we  obtain  the following

\begin{cor}
\label{cor:compatible}
Let $r: \frakh^* \lon \wedge^2 \frakg$ be a 
nondegenerate triangular dynamical $r$-matrix,
and $M=\frakh^* \times G$ its associated symplectic manifold.
Let $\nabla$ be the symplectic connection on $M$ as in
Theorem \ref{thm:connection}. 
Suppose that  $\Omega =\omega +\hbar \omega_{1} +\cdots +\hbar^{i}\omega_{i}
+\cdots \in Z^{2}(M)^{G} \flb\hbar \frb$
 satisfies the condition that $i_{\vh{}}\omega_{i}=0,  \ \forall i\geq 1, \ h\in \frakh$.
Then the Fedosov star product on $M$ corresponding to $(\nabla , \Omega )$
is a  compatible star product.
\end{cor}

 Combining with Theorem \ref{thm:star},
we are  lead to the  following main result of the paper.

\begin{thm}
Any nondegenerate triangular dynamical $r$-matrix is
quantizable.
\end{thm}

More generally, if $r$ is a splittable triangular dynamical
$r$-matrix, according to Corollary \ref{cor:conn-Poisson},
 the corresponding Poisson manifold
$M=\frakh^* \times G$ admits a $G\times H$-invariant
leafwise (w.r.t. the symplectic foliation) Poisson connection
 such that 
$\nabla_{X}\vh{}=0, \forall h\in \frakh$.
Applying Theorem  \ref{thm:Fedosov} leafwisely, 
we thus have the following

\begin{thm}
Any splittable  triangular dynamical $r$-matrix is
quantizable.
\end{thm}

The rest of  the section is devoted to the proof of Theorem \ref{thm:Fedosov}.
We will start with the following

\begin{pro}
\label{pro:propertity}
Under the same hypothesis as in Theorem \ref{thm:Fedosov}, 
we have
\begin{enumerate}
\item For any  $(r, s)$-type tensor 
 $S\in  \calt^{(r, s)}M$ and   $ h\in \frakh$, we have
$\nabla_{\vh{}}S=L_{\vh{}}S$.
\item $\forall X\in {\frak X}(M)$ and $i \geq 1$,
 $\nabla_{X}(J^*d\lambda^{i} ) =0$. 
\item Given any  $\theta \in \Omega^{1} (M)$, if $\vh{}\per \theta =0$,
then $\vh{}\per \nabla_{X} \theta =0$, $\forall X\in {\frak X} (M)$.
\item $R(X, Y, Z, W)=0$, if any    of the vectors $X, Y, Z, W$ is
tangent to the $H$-orbits.
\item $\nabla_{\vh{}}R=0$,  $\forall h\in \frakh$.
\end{enumerate}
\end{pro}
\pf (i). Since $\nabla$ is torsion-free, for any vector field
$X\in {\frak X} (M)$, we have  
$$\nabla_{\vh{}} X=\nabla_X {\vh{}} +[\vh{}, \ X]=
[\vh{}, \ X]= L_{\vh{}}X. $$
 This implies that $\nabla_{\vh{}} \theta
=L_{\vh{}}\theta$ for any one form $\theta \in \Omega^{1}(M)$.
Therefore $\nabla_{\vh{}}S=L_{\vh{}}S$ for any $(r, s)$-type
tensor $S\in \calt^{(r, s)}M$.

(ii). Since $J: M \lon \frakh^*$ is a momentum map, it follows that
$J^* d\lambda^{i}  =\omega^{b} \vh{i}$, where
$\omega^{b}: {\frak X}  (M)\lon \Omega^{1}(M)$ is
the isomorphism induced by the symplectic
structure $\omega$. Hence
$\nabla_{X}(J^*d\lambda^{i})
=\nabla_{X} (\omega^{b} \vh{i})=\omega^{b} (\nabla_{X}\vh{i})=0$, 
since $\nabla$ is a symplectic connection.

(iii). We have  $\nabla_{X}(\vh{}\per \theta )=(\nabla_{X}\vh{})\per \theta
+\vh{} \per \nabla_{X}\theta =\vh{} \per \nabla_{X}\theta$.
The claim thus follows.

(iv). Let $\Phi$ denote the  $H$-action on $M$. For any
$h\in \frakh$, since $\nabla$ is $H$-invariant, it follows that
$\forall W, Y\in {\frak X}  (M)$,
$ \nabla_{(\Phi_{\exp{th}*}W )}(\Phi_{\exp{th}*}Y ) =
\Phi_{\exp{th}*}(\nabla_W Y)$. Taking the derivative at $t=0$,
one obtains that
$$\nabla_{[\vh{}, W]}Y+\nabla_{W}[\vh{}, Y]=[\vh{}, \ \nabla_{W}Y]. $$
Hence, 
\be
R(\vh{}, W)Y&=&\nabla_{\vh{}} \nabla_{W}Y-\nabla_{W}\nabla_{\vh{}} Y-
\nabla_{[\vh{}, W]}Y\\
&=&[\vh{}, \ \nabla_{W}Y] -\nabla_{W}[\vh{}, Y]-\nabla_{[\vh{}, W]}Y\\
&=&0.
\ee

On the other  hand, we know that  $R(Z, W)\vh{} =0$,
 since $\vh{}$ is parallel by
assumption. This means that $R(X, Y, Z, W)=0$ if $Y=\vh{}$  or
$Z=\vh{}$. Since $R(X, Y, Z, W)$ is antisymmetric with respect to
$W, Z$, and symmetric with respect to $X, Y$
according to Proposition \ref{pro:curvature},
 the conclusion thus follows. 

(v). Since both  the connection $\nabla$ and the symplectic structure $\omega$
are $H$-invariant, the symplectic curvature $R$, as defined
by Equation (\ref{eq:sym-curvature}), 
is also $H$-invariant. Hence, for any $h\in \frakh$,
according to (i), $\nabla_{\vh{}}R=L_{\vh{}}R=0$.

This completes the proof of the proposition. \qed

By $K\subset TM$, we denote the  integrable distribution  on $M$
corresponding to the $H$-orbits, and $K^{\perr}$ its
conormal subbundle.  That is, a covector $\theta$ is in $K^{\perr}$
iff $<\theta , \vh{}>=0, \ \forall h\in \frakh$. For
any $x\in M$,   by
$\mbox{pol} (K_{x}^{\perr} )$, we denote the polynomials  on  $T_{x}M$
generated  by those linear functions corresponding to  covectors
in $K_{x}^{\perr} $.
 By $W_{x}^{\perr}$, we denote the formal power series
in $\hbar$ with coefficients in $\mbox{pol}(K_{x}^{\perr} )$. Clearly
$W_{x}^{\perr}$ is a subalgebra  of the  Weyl algebra $W_{x}$.
Let $W^{\perr}=\cup_{x\in M}W_{x}^{\perr}$ be the subbundle
of $W$. We  also consider   $W^{\perr}\ot \wedge^q 
K^{\perr}$, a subbundle of $W\ot \wedge^q T^* M$, whose
space of  sections  is denoted by $\gm W^{\perr}\ot (\Lambda ^{\perr})^q$.
As before, let us  fix a basis $\{h_{1}, \cdots , h_{l}\}$ of
$\frakh$, and denote by $(\lambda^{1}, \cdots , \lambda^l )$
its induced coordinate system on $\frakh^*$.  Since
$J: M\lon \frakh^* $ is a momentum map,  we have
$X_{J^*\lambda^i}=\vh{i}, \ \forall
i=1, \cdots, l$. It thus follows that $J_{*}\vh{i}= J_{*} X_{J^*\lambda^i}=
0, \ \forall i=1, \cdots, l$,  since $\frakh$ is Abelian.
Next we need to  extend $\{\vh{1}, \cdots ,\vh{l}\}$
to a set of (local) vector fields which constitutes a basis of tangent 
fibers of $M$.
  For this purpose, let  $\{u_{1}, \cdots , u_{m}\}$ be   (local)
vector fields on $M$ 
 tangent to the J-fibers such that $\{\vh{1}, \cdots ,\vh{l},
u_{1}, \cdots , u_{m}\}$ constitutes  a basis of the tangent spaces of the J-fibers.
Choose (local) vector fields  $\{v_{1}, \cdots , v_{l}\}$  on $M$ such that
$J_{*}v_{i}=\parr{}, \ \forall i=1, \cdots, l$, which is
always possible since $J$  is a submersion. It is
easy to see that locally  $\{\vh{1}, \cdots ,\vh{l}, v_{1}, \cdots , v_{l},
u_{1}, \cdots , u_{m}\}$ constitutes  a basis of the  tangent 
fibers of $M$.
Let  $ \{ \vhs{1}, \cdots ,\vhs{l}, \vs{1}, \cdots \vs{l},
\us{1}, \cdots , \us{m}\}$ be  its dual basis. 
Then any section of $W\ot \Lambda $ can be written as 
\begin{equation}
a=\sum \hbar^{k}   a_{k, i_{1}\cdots i_{p}, j_{1}\cdots j_{q}}
y^{i_{1}}_{*}\cdots y^{i_{p}}_{*}x^{j_{1}}_{*} \wedge \cdots \wedge
 x^{j_{q}}_{*},
\end{equation}
where all $y^{i}_{*}$'s and $x_{*}^{i}$'s are  either
$\vhs{i}, \vs{i}$ or $\us{i}$, and
 the coefficients $a_{k, i_{1}\cdots i_{p}, j_{1}\cdots j_{q}} $
are   covariant  tensors symmetric with respect to
$i_{1} \cdots i_{p}$ and antisymmetric in $j_{1}\cdots j_{q}$.
It is simple to see that 
a section  $a$ belongs to $ \gm W^{\perr}\ot (\Lambda ^{\perr})^{q}$
 iff there are no terms involving explicit $\hs{i}$'s   in the above expression.

\begin{lem}
\label{lem:basic-relation}
\begin{enumerate}
\item For any $i=1, \cdots , l$, $J^{*}d \lambda^{i}=\vs{i}$;
\item for any $i, j$, $ \nabla_{\vh{i}} \vs{j}=0$,
and $\nabla_{\vh{i}}\hs{j}$ and  $\nabla_{\vh{i}} \us{j}$ 
belong to $\gm K^{\perr}$;
\item for any $i, j$, $ \pi (\vs{i} ,\hs{j}) =\delta_{ij}, \ \pi (\vs{i} ,\vs{j})=0, \ \ 
\pi (\vs{i} ,\us{j})=0$;
\item the commutatant of $\{\vs{1}, \cdots , \vs{l}\}$ in $\gm W$
is $\gm W^{\perr}$.
\end{enumerate}
\end{lem}
\pf (i) $<J^{*}d \lambda^{i}, \ v_{j}>=<d \lambda^{i},   \ J_{*} v_{j}>
=<d \lambda^{j},  \ \ \parr{}>=\delta_{ij}$. Similarly, we have
$<J^{*}d \lambda^{i}, u_{j}>=0$ and $<J^{*}d \lambda^{i}, h_{j}>=0$.
Therefore, $J^{*}d \lambda^{i}=\vs{i}$.

(ii) According to Proposition \ref{pro:propertity}, 
 $\nabla_{\vh{i}} \vs{j}=\nabla_{\vh{i}}(J^* d\lambda^{j})=0$.
Also, $\forall k$, $<\nabla_{\vh{i}}\hs{j}, \vh{k}>$\\
$= \nabla_{\vh{i}}<\hs{j}, \vh{k}>-<\hs{j}, \nabla_{\vh{i}} \vh{k}>=0$.
Hence it follows that $\nabla_{\vh{i}}\hs{j}\in \gm   K^{\perr}$.
Similarly, we can  prove that $\nabla_{\vh{i}} \us{j} \in \gm K^{\perr}$.

(iii) We have  $ \pi (\vs{i} ,\hs{j}) =  
<\pi^{\#}{(J^{*}d \lambda^{i})}, \ \hs{j}>
=<\vh{i},  \ \hs{j}>=\delta_{ij}$. Similarly, we can show that 
$\pi (\vs{i} ,\vs{j})=0$ and $ \pi (\vs{i} ,\us{j})=0$.

(iv) Assume that $a\in \gm W$ such that
 $[a, \ \vs{i}]=0, \ \forall i=1, \cdots ,l$.
It thus follows that $\{a, \ \vs{i}\}=0$, where the Poisson
bracket refers to the one corresponding to the
fiberwise symplectic structure on $TM$. Thus $a\in \gm W^{\perr}$ according
to (iii). \qed

\begin{lem}
\label{lem:4.3}
\begin{enumerate}
\item  $\gm W^{\perr}\ot\Lambda ^{\perr}$ is closed under the multiplication
$\circ$ as defined by  Equation (\ref{eq:circle}).
\item  $\gm W^{\perr}\ot\Lambda ^{\perr}$ is closed under both the
operators $\delta$ and $\delta^{-1}$.
\item $\gm W^{\perr}\ot\Lambda ^{\perr}$ is invariant under the covariant
derivative $\nabla_{X}, \ \forall X\in {\frak X}  (M)$.
\end{enumerate}
\end{lem}
\pf (i) and (ii) are obvious. For (iii), note that $\gm (K^{\perr})$ is
invariant under the covariant derivative $\nabla_{X} $ according to
Proposition \ref{pro:propertity} (iii). Hence
$\gm W^{\perr}\ot\Lambda ^{\perr}$ is also invariant. \qed

As an immediate consequense, we have the  following

\begin{cor}
\label{cor:par}
If $a\in \gm W^{\perr}\ot\Lambda ^{\perr}$ and $\nabla_{\vh{}}a=0, \forall
h\in \frakh$, then $\partial a\in \gm W^{\perr}\ot\Lambda ^{\perr}$.
\end{cor}

To prove  Theorem \ref{thm:Fedosov},  we start with the following

\begin{lem}
\label{lem:ga0}
Under the same hypthesis as in Theorem \ref{thm:Fedosov},
 we have $ \gamma_{0}=
\delta^{-1}\tilde{\Omega}\in \gm W^{\perr}\ot\Lambda ^{\perr}$ and
$\nabla_{\vh{}}\gamma_{0}=0, \forall h\in \frakh$.
\end{lem}
\pf According to Equation (\ref{eq:tildeOme}), 
 we know that $\tilde{\Omega}= \Omega -\omega +R=
R+ \hbar \omega_{1}+\hbar^{2}\omega_{2}+\cdots $. By assumption,
 we have  $\omega_{i}\in \gm W^{\perr}\ot\Lambda ^{\perr}, \ \forall
i\geq 1 $.  On the other hand, according to Proposition
\ref{pro:propertity} (iv), we know that
$R\in \gm W^{\perr}\ot\Lambda ^{\perr}$. Therefore,
$\tilde{\Omega}\in \gm W^{\perr}\ot\Lambda ^{\perr}$. Hence
$\gamma_{0}\in \gm W^{\perr}\ot\Lambda ^{\perr}$ by Lemma \ref{lem:4.3}.

Finally, note that  for any
$h\in \frakh$, $L_{\vh{}} \omega_{i} = i_{\vh{}} (d\omega_{i})+d(i_{\vh{}}
\omega_{i})=
d(i_{\vh{}}\omega_{i} )=0$. According to Proposition
\ref{pro:propertity}, we have  $L_{\vh{}}R
=\nabla_{\vh{}}R=0$. Hence $L_{\vh{}}\tilde{\Omega}=0$.
It thus follows that $\nabla_{\vh{}}\gamma_{0}=L_{\vh{}}\gamma_{0}
=L_{\vh{}} \delta^{-1}\tilde{\Omega}=\delta^{-1}L_{\vh{}} \tilde{\Omega}=0$.\qed

\begin{pro}
\label{pro:gamma}
Under the same hypthesis as in Theorem \ref{thm:Fedosov}, 
 the element  $ \gamma$, defined as in Theorem
\ref{label:f1}, belongs to $ \gm W^{\perr}\ot\Lambda ^{\perr}$ and 
satisfies
$\nabla_{\vh{}}\gamma =0, \forall h\in \frakh$.
\end{pro}
\pf  We prove this proposition  by induction. Assume that
 $\gamma_n \in \gm W^{\perr}\ot\Lambda ^{\perr}$ 
and $\nabla_{\vh{}}\gamma_{n} =0, \forall h\in \frakh$. It suffices
to show that $\gamma_{n+1}$   satisfies the same conditions.
By Equation (\ref{eq:r-iteration2}), $\gamma_{n+1}$ and $\gamma_n$
 are related by  the following equation: 
\begin{equation}
\gamma_{n+1} =\gamma_{0}+\delta^{-1}(\partial \gamma_{n} +\ih \gamma_{n}^{2} ),
\ \ \ \ \forall n\geq 0.
\end{equation}

According to Corollary \ref{cor:par}, we have  $\partial \gamma_n \in 
\gm W^{\perr}\ot\Lambda ^{\perr}$. On the other hand,
 by Lemma \ref{lem:4.3},
$\gamma_{n}^{2} \in \gm W^{\perr}\ot\Lambda ^{\perr}$.
Hence $\gamma_{n+1}\in \gm W^{\perr}\ot\Lambda ^{\perr}$ according
to  Lemma \ref{lem:4.3}  and Lemma \ref{lem:ga0}.

Now
\be
&& \nabla_{\vh{}}(\gamma_{n+1} )\\
&=& \nabla_{\vh{}}\gamma_{0}+\nabla_{\vh{}}\delta^{-1}(\partial \gamma_{n} +\ih \gamma_{n}^{2} )\\
&=& L_{\vh{}}\delta^{-1}(\partial \gamma_{n} +\ih 
\gamma_{n}^{2} )\\
&=& \delta^{-1} ( L_{\vh{}} \partial \gamma_{n} +\ih L_{\vh{}} \gamma_{n}^{2} )\\
&=& 0.
\ee
Here, in the last step,
 we used the relation $ L_{\vh{}} \partial =\partial  L_{\vh{}}$, which
follows from the fact that the symplectic connection
is $H$-invariant. This concludes the proof.  \qed

As in Appendix, for any $a\in C^{\infty}(M)$, we denote by
$\tilde{a}\in W_{D}$ its parallel lift, i.e., $D\tilde{a}=0$ and
$\tilde{a}|_{y=0}=a$. Theorem \ref{thm:Fedosov}
is in fact an immediate consequence of the following
\begin{pro}
\label{pro:Fed}
Under the same hypothesis as in Theorem \ref{thm:Fedosov}, 
\begin{enumerate}	
\item if $a=J^* f$ for  $f\in C^{\infty}(\frakh^* )$,
then 
\begin{equation}
\label{eq:a}
\tilde{a}=\sum_{k=0}^{\infty}\frac{1}{k!}
J^* (\frac{\partial^{k} f}{\partial \lambda^{i_{1}}\cdots \partial 
\lambda^{i_{k}} } )\vs{i_{1}}\cdots \vs{i_{k}};
\end{equation}
\item for any $a\in C^{\infty}(M)$,
$$\tilde{a}=\sum_{k=0}^{\infty}\frac{1}{k!} 
(\vh{i_{1}}\cdots \vh{i_{k}}  a )\hs{i_{1}}\cdots \hs{i_{k}} +T,$$
where the reminder $T$   does  not contain any terms
which are   pure polynomials of $\hs{i}$'s.
\end{enumerate}
\end{pro}
\pf For (i), it suffices to prove that $\tilde{a}$ given
by Equation (\ref{eq:a}) is a parallel section.  According
 to Proposition \ref{pro:gamma} and Lemma \ref{lem:basic-relation},
we have   $[\gamma , \tilde{a} ]=0$. Thus it follows that
$D\tilde{a}=-\delta \tilde{a}+\partial \tilde{a}$,
which  clearly vanishes  since $\partial \vs{i}=0$ by 
Proposition \ref{pro:propertity} (ii) and Lemma \ref{lem:basic-relation} (ii).

For (ii), recall that $\tilde{a}$ is determined  by the iteration
formula  
\begin{equation}
\label{eq:a-iteration}
a_{n+1} =a+\delta^{-1}(\partial {a_{n}} +[\ih \gamma ,
{a_{n}} ]).
\end{equation}

So it suffices to prove  that 
$$a_{n}=\sum_{k=0}^{n}
\frac{1}{k!} (\vh{i_{1}}\cdots \vh{i_{k}}a_{0} )\hs{i_{1}}\cdots \hs{i_{k}} 
+T_{n},$$
where  each term in the reminder $T_{n}$   is not a pure polynomial
of $\hs{i}$'s.  This can be proved by induction again.

Assume that this assertion holds for $a_{n}$. To show that it still
 holds
for $a_{n+1}$,  we need to analyze which terms in $a_{n}$
would produce pure polynomials of $\hs{i}$'s out of
Equation (\ref{eq:a-iteration}). Since
$\gamma \in  \gm W^{\perr}\ot\Lambda ^{\perr}$, we may ignore
$\delta^{-1} [\ih \gamma , {a_{n}} ] $ and  only consider
  $\delta^{-1} \partial {a_{n}}= \delta^{-1} (
 \sum_{i}\nabla_{\vh{i}}a_{n} \wedge \hs{i}
+ \sum_{i}\nabla_{v_{i}}a_{n} \wedge \vs{i}
+\sum_{i}\nabla_{u_{i}}a_{n} \wedge \us{i})$.
 From this, it is clear that those terms containing
 pure polynomials of $\hs{i}$'s arise
only from $ \delta^{-1} (\sum_{i}\nabla_{\vh{i}}a_{n} \wedge
 \hs{i})$.
Now a general term in $a_n$ has the form $ \hbar^{k}
a_{\alpha \beta \gamma } (x) \vs{\alpha}\hs{\beta}\us{\gamma }$,
where $\alpha , \beta$ and $\gamma $ are multi-indexes.  However,
\be
&&\nabla_{\vh{i}} (a_{\alpha \beta \gamma } (x) 
\vs{\alpha}\hs{\beta}\us{\gamma })\\
&=&(\vh{i}a_{\alpha \beta \gamma } (x) )\vs{\alpha}\hs{\beta}\us{\gamma }+
a_{\alpha \beta \gamma }(x)(\nabla_{\vh{i}}\vs{\alpha}) \hs{\beta}\us{\gamma }
+a_{\alpha \beta \gamma }(x)\vs{\alpha} (\nabla_{\vh{i}} \hs{\beta})\us{\gamma }
+a_{\alpha \beta \gamma }(x)\vs{\alpha}  \hs{\beta} (\nabla_{\vh{i}}\us{\gamma })\\
&=&(\vh{i}a_{\alpha \beta \gamma } (x))\vs{\alpha}\hs{\beta}\us{\gamma }
+a_{\alpha \beta \gamma }(x) \vs{\alpha} (\nabla_{\vh{i}} \hs{\beta})\us{\gamma }
+a_{\alpha \beta \gamma } (x)\vs{\alpha}  \hs{\beta} (\nabla_{\vh{i}}\us{\gamma }).
\ee
According to Lemma \ref{lem:basic-relation},
 neither   $\nabla_{\vh{i}} \hs{\beta}$ nor
$\nabla_{\vh{i}} \us{\gamma }$ will be a pure polynomial
 of $\hs{i}$'s. Hence to produce a pure $\hs{i}$-polynomial
term, one needs that $\alpha=\gamma =0$. And in this case,
the resulting pure $\hs{i}$-polynomial
term is $\hbar^{k} (\vh{i}a_{0 \beta 0}(x) )\hs{\beta}$. In conclusion, only
pure $\hs{i}$-polynomial terms in $a_n$ can give rise to
pure $\hs{i}$-polynomial terms in   $\delta^{-1} \partial {a_{n}}$.
Hence the pure $\hs{i}$-polynomial terms in $a_{n+1}$ is
$a_{0}+\sum_{k=0}^{n} \frac{1}{k!}\frac{1}{k+1}\vh{i}
 (\vh{i_{1}}\cdots \vh{i_{k}} a_{0})\hs{i} \hs{i_{1}}\cdots \hs{i_{k}} $,
which clearly  equals to $\sum_{k=0}^{n+1}
\frac{1}{k!} (\vh{i_{1}}\cdots \vh{i_{k}}a_{0} )\hs{i_{1}}\cdots \hs{i_{k}} $.
This concludes the proof.
\qed

\section{Classification}
This section is devoted to the classification of quantization 
of a non-degenerate triangular dynamical $r$-matrix.
Our method relies heavily on the classification result 
of star products on a symplectic manifold.
First, let us introduce the following:

\begin{defi}
\label{defi:equi}
Two quantizations $F(\lambda )$ and $E(\lambda )$ of a triangular 
dynamical $r$-matrix are said to be equivalent if there exists a
 $T(\lambda ) : \frakh^* \lon  (U\frakg )^{\frakh} \flb\hbar \frb$
 satisfying  the condition that $T(\lambda )\eqq 1 (\mbox{mod } \hbar)$
and $\epsilon ( T(\lambda ) )=1$ such that
\begin{equation}
\label{eq:G}
E(\lambda )=\Delta T(\lambda)^{-1} 
F(\lambda )T_{1}(\lambda -\half \hbar h^{(2)})
T_{2}(\lambda +\half \hbar h^{(1)}).
\end{equation}
\end{defi}

To justify this   definition, we need the following result,
which interprets this equivalence in terms of star products.

\begin{thm}
\label{thm:transformation}
Given a compatible star
product $*_{\hbar}$ on the Poisson manifold $(M, \pi )$
associated to a triangular dynamical $r$-matrix  $r(\lambda )$, assume that
$T(\lambda ) : \frakh^* \lon  (U\frakg )^{\frakh} \flb\hbar \frb$
 satisfies  the condition that $T(\lambda )\eqq 1 (\mbox{mod } \hbar)$
and $\epsilon  (T(\lambda ) )=1$. Then the $*$-product:
\begin{equation}
\label{eq:T}
f\tilde{*_{\hbar}}g=\Vec{T}^{-1}(\Vec{T}f *_{\hbar} \Vec{T}g ), \ \ \ 
\forall f , g\in C^{\infty}(M)
\end{equation}
is still a compatible star-product. Moreover,
 if $f, g\in C^{\infty}(G)$,
\begin{equation}
f\tilde{*_{\hbar}}g=\Vec{E(\lambda )}(f, g),
\end{equation}
where $E(\lambda )$ is given by Equation (\ref{eq:G}).
\end{thm}

Thus we are lead to the following

\begin{defi}
Compatible star-products $*_{\hbar}$ and $\tilde{*_{\hbar}}$ are
said to be strongly equivalent iff they are related
by Equation (\ref{eq:T}) for
 some $T(\lambda ): \frakh^* \lon  (U\frakg )^{\frakh} \flb\hbar \frb$
satisfying the property that $T(\lambda )\eqq 1 (\mbox{mod } \hbar)$
and $\epsilon  (T(\lambda ) )=1$.
\end{defi}

An immediate consequence of Theorem \ref{thm:transformation} is
the following:

\begin{cor}
If $F(\lambda )$ is a  quantization  of a triangular
dynamical $r$-matrix $r: \frakh^* \lon \wedge^{2}\frakg$
 and  $T(\lambda ) : \frakh^* \lon  (U\frakg )^{\frakh} \flb\hbar \frb$ 
satisfies  the condition that 
$T(\lambda ) \eqq 1 (\mbox{mod } \hbar)$ and $\epsilon  (T(\lambda ) )=1$,
 then $$E(\lambda )=\Delta T(\lambda)^{-1}
F(\lambda )T_{1}(\lambda -\half \hbar h^{(2)})
T_{2}(\lambda +\half \hbar h^{(1)})$$
is also a  quantization of $r(\lambda )$.
\end{cor}

Due to  this fact, Definition  (\ref{defi:equi}) is well justified.  Indeed,
Theorem \ref{thm:transformation} allows us to reduce
the classification problem of quantizations of
a triangular dynamical $r$-matrix to that of
strongly equivalent star products on $M$.\\\\
{\bf Remark} Let $R_{E}(\lambda )=E^{21}(\lambda )^{-1}
E^{12}(\lambda )$ and $R_{F}(\lambda )=F^{21}(\lambda )^{-1}
F^{12}(\lambda )$. It is easy to see that they are
related by
\begin{equation}
\label{eq:R}
R_{E}(\lambda )=T_{2}(\lambda -\half \hbar h^{(1)})^{-1}
T_{1}(\lambda +\half \hbar h^{(2)})^{-1}R_{F}(\lambda )
T_{1}(\lambda -\half \hbar h^{(2)})
T_{2}(\lambda +\half \hbar h^{(1)}).
\end{equation}
Alternatively, we may define a quantization of a triangular
dynamical $r$-matrix $r(\lambda )$ to be
an element $R(\lambda )  =1+\hbar r(\lambda )+\cdots   
 \in U \frakg \ot U \frakg  \flb\hbar \frb$
satisfying  the QDYBE,
and  define an equivalence of  quantizations
 by Equation (\ref{eq:R}).
This definition sounds  weaker than our original one.
We  do not know, however, at this moment
 whether these two definitions are equivalent.
 It would be interesting to have this clarified.

To prove Theorem \ref{thm:transformation}, we need a lemma.

\begin{lem}
\label{lem:6.4}
Assume that $T(\lambda ) :
 \frakh^* \lon  (U\frakg)^{\frakh} \flb\hbar \frb$ is
as in Theorem \ref{thm:transformation}, then
\begin{enumerate}
\item  $(ad_{\theta})^{n}(T\ot 1 )=(-\half )^{n}
\sum_{i_{1}\cdots i_{n}} \frac{\partial^{n} T}{\partial \lambda^{i_{1}}
\cdots \partial \lambda^{i_{n}}} \ot h_{i_{1}}\cdots h_{i_{n}} $;
\item $\Theta (T\ot 1 ) \Theta^{-1} = T_{1}(\lambda -\half \hbar h^{(2)})$;
\item $\Theta (1\ot T)\Theta^{-1} = T_{2}(\lambda +\half \hbar h^{(1)})$;
\item $\Theta (T\ot T) \Theta^{-1}
=T_{1}(\lambda -\half \hbar h^{(2)})T_{2}(\lambda +\half \hbar h^{(1)})$.
\end{enumerate}
\end{lem}
\pf (i) We prove this equation  by induction. Obviously, it holds for
$n=0$. Assume that it holds for $n=k$. Now
\be
&&(ad_{\theta})^{k+1}(T\ot 1 )\\
&=& ad_{\theta}[ (-\half )^{k}
\sum_{i_{1}\cdots i_{k}} \frac{\partial^{k} T}{\partial \lambda^{i_{1}}
\cdots \partial \lambda^{i_{k}}} \ot h_{i_{1}}\cdots h_{i_{k}} ]\\
&=&(-\half )^{k} \half \sum_{i_{1}\cdots i_{k}} \sum_{i} 
([h_{i}\ot \parr{},  \ 
\frac{\partial^{k} T}{\partial \lambda^{i_{1}}
\cdots \partial \lambda^{i_{k}}} \ot h_{i_{1}}\cdots h_{i_{k}} ]\\
&&-[ \parr{} \ot h_{i}, \ \frac{\partial^{k} T}{\partial \lambda^{i_{1}}
\cdots \partial \lambda^{i_{k}}} \ot h_{i_{1}}\cdots h_{i_{k}} ])\\
&=&(-\half )^{k+1} \sum_{i_{1}\cdots i_{k+1}} \frac{\partial^{k+1} T}{\partial 
\lambda^{i_{1}} \cdots \partial \lambda^{i_{k+1}}} \ot h_{i_{1}}\cdots 
h_{i_{k+1 }}. 
\ee

(ii)
We have 
\be 
&&\Theta (T\ot 1 ) \Theta^{-1}\\
&=&\exp{(\hbar ad_{\theta} )} (T\ot 1 )\\
&=&\sum_{k=0}^{\infty} \frac{1}{k!}(\hbar ad_{\theta})^{k}(T\ot 1 )\\
&=&\sum_{k=0}^{\infty}\frac{1}{k!}(-\frac{\hbar}{2})^{k}
\frac{\partial^{k} T}{\partial \lambda^{i_{1}}
\cdots \partial \lambda^{i_{k}}} \ot h_{i_{1}}\cdots h_{i_{k}} \\
&=& T_{1}(\lambda -\half \hbar h^{(2)}).
\ee

(iii) is proved  similarly,  and (iv) follows from (ii) and (iii). \qed
{\bf Proof of Theorem \ref{thm:transformation}}
 If $f, g\in C^{\infty}(\frakh^* )$,
then $\Vec{T}f=f$ and $\Vec{T}g=g$ since $\epsilon (T)=1$. 
Hence 
$$\Vec{T}^{-1}(\Vec{T}f*_{\hbar} \Vec{T}g)=fg .$$ 

Now if $f\in C^{\infty}(\frakh^* )$ and $g\in C^{\infty}(G)$,
\be
&& \Vec{T}f*_{\hbar} \Vec{T}g\\
&=&f*_{\hbar} \Vec{T}g\\
&=&\Vec{\Theta} (f, \Vec{T}g)\\
&=&\Vec{\Theta (1\ot T)}(f, g) \ \ \mbox{ (by Lemma \ref{lem:6.4})}\\
&=&\Vec{ T_{2}(\lambda +\half \hbar h^{(1)})\Theta }(f, g)\\
&=&\sum_{k=0}^{\infty} 
 \frac{1}{k!}(\frac{\hbar}{2})^{k}  \Vec{(h_{i_{1}}\cdots h_{i_{k}} \ot
\frac{\partial^{k} T}{\partial \lambda^{i_{1}}
\cdots \partial \lambda^{i_{k}}} )\Theta }(f, g)\\
&=&\Vec{(1\ot T)\Theta }(f, g)\\
&=&\Vec{T}(f*_{\hbar}g).
\ee
Here  in the last equality, we used the fact that $\Vec{T}$
does not involve any  derivative $\parr{}$.
 So we have proved that
   $\Vec{T}^{-1}(\Vec{T}f*_{\hbar} \Vec{T}g)
=\Vec{\Theta }(f, g)$.

Finally, assume that $f, g\in C^{\infty}(G)$. According to Theorem
\ref{thm:full-star}, 
\be
&&\Vec{T}f*_{\hbar} \Vec{T}g\\
&=&\Vec{(F(\lambda )\Theta })(\Vec{T}f, \Vec{T}g )\\
&=&\Vec{F(\lambda )\Theta (T\ot T)} (f, g) \ \ \mbox{ (by Lemma \ref{lem:6.4})}\\
&=&\Vec{F(\lambda )T_{1}(\lambda -\half \hbar h^{(2)})
T_{2}(\lambda +\half \hbar h^{(1)}) \Theta}(f, g).
\ee
It thus follows that 
\be
&&\Vec{T}^{-1}(\Vec{T}f*_{\hbar} \Vec{T}g)\\
&=&\Vec{\Delta T(\lambda )^{-1} F(\lambda )T_{1}(\lambda -\half \hbar h^{(2)})
T_{2}(\lambda +\half \hbar h^{(1)}) \Theta}(f, g)\\
&=&\Vec{E(\lambda )  \Theta}(f, g).
\ee
This concludes the proof. \qed.

The rest of   the  section is devoted to the  classification
of strongly equivalent classes  of
compatible star products on $M=\frakh^* \times G$.
 The classification
of star products on  a  general symplectic manifold was studied
by many authors,  for example, see 
\cite{BCG, D, NT1, WeinsteinX:1996, Xu0}. Here 
we follow the  elementary approach due to
 Bertelson, Bieliavsky and Gutt \cite{BBG} concerning invariant star products.

First we prove 

\begin{thm}
\label{thm:Fed-classification}
Let $M=\frakh^* \times G$ be the symplectic
manifold corresponding  to a non-degenerate dynamical
$r$-matrix $r(\lambda )$. Two compatible Fedosov $*$-products
are strongly  equivalent iff their Weyl curvatures
$\Omega_{*}$ and $\Omega$ are strongly  cohomologous,
i.e., $\Omega_{*}-\Omega=d\theta $, where $\theta \in \Omega^{1}(M)\flb\hbar \frb$
is $G\times H$-invariant and
satisfies  $i_{\vh{}}\theta =0$, $\forall h\in \frakh$.
\end{thm}

From now on, in this section, by $M$ we always mean 
the symplectic manifold $\frakh^* \times G$ associated
with a non-degenerate dynamical $r$-matrix.
Let $\frakg =\frakh \oplus \muu $ be a reductive decomposition
as in  Lemma \ref{lem:decom}, and  $\{h_{1}, \cdots , h_{l}\}$
a basis in $\frakh$, and  $\{e_{1},  \cdots , e_{m}\}$ 
a basis of $\muu$. If we choose $v_{i}=\frac{\partial}{\partial \lambda^i}$
and $u_{i}=\ve{i}$, then $\{\vh{1},  \cdots , \vh{l}, v_{1}, \cdots ,
v_{l}, u_{1}\cdots , u_{m}\}$ constitutes  a local (in fact global in
this case) basis of tangent fibers of $M$, which satisfies all the
required properties as in the construction preceding
Lemma \ref{lem:basic-relation}. In what follows, we will fix such a
choice, and  denote by  $ \{ \vhs{1}, \cdots ,\vhs{l}, \vs{1}, \cdots \vs{l},
\us{1}, \cdots , \us{m}\}$ its dual basis.

\begin{lem}
\label{lem:6.7}
Assume  that $D$ is an Abelian connection
defining a compatible $*$-product on $M$
as in Corollary \ref{cor:compatible}.
For any $a\in C^{\infty}(M)$, let 
\begin{equation}
\label{eq:tilde_a}
\tilde{a}=\sum \hbar^{k}D_{k, \alpha  \beta  \gamma }(a) \vs{\alpha}
\hs{\beta}\us{\gamma} \in \gm (W)
\end{equation}
be its parallel lift, where $\alpha, \beta$ and
$\gamma$ are multi-indexes, and  $D_{k, \alpha  \beta  \gamma }$
are certain  differential operators on $M$.
If an  operator $D_{k, \alpha  \beta  \gamma }$ 
 involves a derivative of  $\lambda \in \frakh^* $, then  the corresponding
term  satisfies $|\alpha | >0$.
\end{lem}
\pf As  it  is known,  $\tilde{a}$ is given  by the iteration formula
$$
{a_{n+1}}=a_{0}+\delta^{-1}(\partial {a_{n}}+[\ih \gamma , {a_{n}}]), $$
so it suffices to show that ${a_{n}}$ possesses such a property
for any $n$, which we shall prove by induction.

 Assume that  ${a_{n}}$ possesses this property, and we
 need to show that
 so does ${a_{n+1}}$. Let 
$\hbar^{k}D_{k ,\alpha  \beta  \gamma }(a) \vs{\alpha}
\hs{\beta}\us{\gamma} $ be a    term in $a_{n+1}$,
where  $D_{k ,\alpha  \beta  \gamma }$ involves a 
derivative of  $\lambda$. There are two possible
sources that this term may come from. One is from 
$\delta^{-1}[\ih \gamma , {a_{n}}]$. Since this
operation does not affect the part involving derivatives on $a$, so
it must  come from a  term having the form: 
\begin{equation}
\label{eq:res}
\delta^{-1}[\ih \gamma , 
\hbar^{k'}D_{k ,\alpha  \beta  \gamma }(a) \vs{\alpha'}
\hs{\beta'}\us{\gamma'} ],
\end{equation}
 where $\hbar^{k'}D_{k ,\alpha  \beta  \gamma }(a) \vs{\alpha'}
\hs{\beta'}\us{\gamma'}$ is one of the terms in ${a_{n}}$. 
By assumption, we know that $|\alpha' |>0$. 
Since $\gamma  \in \gm W^{\perp}\ot \Lambda^{\perp}$,
it follows from Lemma \ref{lem:basic-relation} that  any  resulting term
in Equation (\ref{eq:res}) has  at least a factor
$\vs{\alpha'}$.

Another possible source is from $\delta^{-1}(\partial {a_{n}})$.
Now 
$$\delta^{-1}\partial {a_{n}}= \sum_{i} (\nabla_{\parr{}} {a_{n}})
\vs{i}+\sum_{i} (\nabla_{\vh{i}} {a_{n}}) \hs{i}
+\sum_{i} (\nabla_{\ve{i}} {a_{n}}) \us{i} .$$
If it arises from the first
term, we are done. Assume that it comes from the second term:
$(\nabla_{\vh{i}} {a_{n}}) \hs{i}$.
Let  $ \hbar^{k}D_{k, \alpha  \beta  \gamma }(a) \vs{\alpha}
\hs{\beta}\us{\gamma}$ be  a  general term in $ {a_{n}}$,
then
$$\nabla_{\vh{i}} (\hbar^{k}D_{k, \alpha  \beta  \gamma }(a) \vs{\alpha}
\hs{\beta}\us{\gamma})=
\hbar^{k}(\vh{i}D_{k, \alpha  \beta  \gamma }) (a) \vs{\alpha}
\hs{\beta}\us{\gamma}
+\hbar^{k}D_{k, \alpha  \beta  \gamma }(a)  \vs{\alpha} (\nabla_{\vh{i}}\hs{\beta})
\us{\gamma} +
\hbar^{k}D_{\alpha , \beta , \gamma }(a)  \vs{\alpha} \hs{\beta}
(\nabla_{\vh{i}}\us{\gamma} ).$$
From this equation, it is clear  that $D_{k, \alpha  \beta  \gamma }$ must
already contain  some derivative  of   $\lambda \in \frakh^*$.
The conclusion  thus follows from the inductive assumption. A similar
argument applies when  it arises from the last term
$(\nabla_{\ve{i}} {a_{n}}) \us{i}$. This concludes the proof. \qed
{\bf Proof of Theorem  \ref{thm:Fed-classification}}  Our
proof here is essentially  a modification of
the proof of  Corollary 5.5.4 in \cite{F2}.

``Necessity." 
 Let
\be
D&=& -\delta +\partial +\ih [\gamma , \cdot ],  \ \ \mbox{and}\\
D_{*}&=& -\delta +\partial +\ih [\gamma_{*}, \cdot ]
\ee
be  the  Abelian connections with Weyl curvatures $\Omega$
and $\Omega_*$, respectively,
 and  $A: W_{D}\lon W_{D_*}$ an  isomorphism of algebras.
 It is standard that $A$ lifts to an automorphism
of the Weyl bundle $W$, which will be denoted by the same
symbol $A: W\lon W$.  Then $A$ is $G\times H$-equivariant.
As in \cite{F2}, we may assume that 
$A (a ) =U\smalcirc a \smalcirc U^{-1}$
 for some $U\in \gm W_{+}$, $\forall a
\in W$.
We may also assume that $U$ is $G\times H$-invariant since
$A$ is $G\times H$-equivariant.
By assumption, we also  know that $A a=a$ if 
$a=\sum_{0}^{\infty}\frac{1}{k!}
\frac{\partial^{k} a_{0}}{\partial 
\lambda^{i_{1}}\cdots \partial \lambda^{i_{k}} }
 \vs{i_{1}}\cdots \vs{i_{k}} \ \forall 
a_{0}\in C^{\infty}(\frakh^* )$,                                       which is the parallel lift of  $a_{0}$
according to Proposition \ref{pro:Fed}.
This implies that $U$ commutes with $\vs{i}, \ i=1, \cdots , l$,
and therefore $U\in \gm W^{\perp}_{+}$ according to 
Lemma \ref{lem:basic-relation}.
Consider another Abelian connection: $\tilde{D}a=(A\smalcirc D\smalcirc A^{-1})
(a)=U\smalcirc D(U^{-1}aU)\smalcirc U^{-1}=
Da-[DU\smalcirc U^{-1}, a]$. Then $\tilde{D}$ has the same Weyl curvature as $D$
(see Theorem 5.5.3 and Corollary 5.5.4 in \cite{F2}), which is
assumed to be $\Omega$. On the other hand,
\begin{eqnarray}
D_{*}a-\tilde{D}a&=&\ih [\gamma_{*}- \gamma -i\hbar (DU\smalcirc U^{-1}), a] \nonumber\\
&=&\ih [\Delta \gamma , a].  \label{eq:difference}
\end{eqnarray}
Since $U\in \gm W^{\perp}_{+}$ and $\gamma_{*}, \gamma \in \gm
W^{\perp}\ot \Lambda^{\perp}$, it follows that $\Delta \gamma\in 
\gm W^{\perp}\ot \Lambda^{\perp}$. It is also
clear that $\Delta \gamma$ is $G\times H$-invariant.
Moreover, from Equation (\ref{eq:difference}), it follows that $[\Delta \gamma , \ a]=0$,
if $a\in W_{D_{*}}$. Hence $\Delta \gamma $ is a scalar form.
Thus $\Omega_{*}-\Omega =d\Delta \gamma$. Clearly, $\Delta \gamma$ is
$G\times H$-invariant and $i_{\vh{}}\Delta \gamma=0, \forall h\in\frakh$. \\\\

``Sufficientity". 
 Assume that $\Omega_* -\Omega =d\theta$, 
$\theta \in \Omega^{1}(M)\flb\hbar \frb$
is $G\times H$-invariant and $i_{\vh{}}\theta =0$, $\forall h\in \frakh$.
Let $\Omega (t)=\Omega +td \theta $ and 
$D_{t}= -\delta +\partial +\ih [\gamma (t) , \cdot ]$ 
be the Abelian connection with Weyl  curvature $\Omega (t)$,
where $\gamma (t)$ is as in Theorem \ref{label:f1} satisfying
$\delta^{-1}\gamma (t)=0$.

Let $H(t) \in \gm W$ be the  solution of  the equation
$ D_{t}H(t)=-\theta +\dot{\gamma}(t) $ satisfying
$H(t)|_{y=0}=0$.
Then $H(t)$ is $G\times H$-invariant since $D_{t}, \ \theta , \ \gamma (t)$
are all $G\times H$-invariant.
On the other hand, since $\gamma (t)\in \gm W^{\perp}\ot \Lambda^{\perp}$ 
according to Proposition \ref{pro:gamma},  and  $\theta \in \gm \Lambda^{\perp}$ by assumption,
it follows that $H(t)\in \gm W^{\perp}$.

According to Theorem 5.5.3 \cite{F2}, the solution
of the Heisenberg equation:
\begin{equation}
\label{eq:Hei}
\frac{d\tilde{a}}{dt}+[H(t), \ \tilde{a}]=0
\end{equation}
establishes an isomorphism $W_{D}\lon W_{D_*}$, which is given  by
$\tilde{a}(0)\lon \tilde{a}(1)$. In fact, $D_{t}\tilde{a}(t)=0$
if $D\tilde{a}(0)=0$.

Clearly, this correspondence is $G\times H$-equivariant since
$H(t)$ is $G\times H$-invariant. So its corresponding
formal differential operator $T: \ C^{\infty}(M) \flb\hbar \frb
\lon C^{\infty}(M) \flb\hbar \frb $ is $G\times H$-invariant.
Finally it remains to show that $T$, 
as a formal differential operator, 
does not involve any derivative of $\lambda \in \frakh^*$.

To show this, for any $a\in C^{\infty}(M)$, let $\tilde{a}\in
W_{D}$ be its parallel lift, and $\tilde{a}(t)$
the solution of Equation (\ref{eq:Hei}) satisfying the initial
condition $\tilde{a}(0)=\tilde{a}$. Then $D_{t}\tilde{a}(t)=0$.
Also, let $a(t)=\tilde{a}(t)|_{y=0}$. Write
$$\tilde{a}(t)=\sum 
\hbar^{k}D_{t, k, \alpha  \beta  \gamma }(a (t))  \vs{\alpha} \hs{\beta}
\us{\gamma}  .$$
If an operator  $ D_{t, k ,\alpha  \beta \gamma }$ involves a 
derivative to $\lambda \in \frakh^*$, we know that $\alpha \neq 0$
according to Lemma \ref{lem:6.7}.
Since $H(t)\in \gm W^{\perp}$, it thus  follows that 
$[H(t), D_{t, k , \alpha  \beta  \gamma }(a (t))  \vs{\alpha} \hs{\beta}
\us{\gamma} ]|_{y=0}=0$.
This implies that $[H(t),  \tilde{a}(t)]|_{y=0}=\DD_{t} a(t)$,
where $\DD_t$ is a  formal differential operator on $M$ involving
no derivatives  of  $\lambda\in \frakh^*$.
Now Equation (\ref{eq:Hei}) implies that
$$\frac{da(t)}{dt}+\DD_{t} (a(t))=0. $$
Therefore the equivalence operator  $T: \ C^{\infty}(M) \flb\hbar \frb
\lon C^{\infty}(M) \flb\hbar \frb $,
 which sends $a(0)$ to $a(1)$, does not involve any
derivative  of $\lambda \in \frakh^*$.
This concludes the proof. \qed

As in \cite{BBG}, by $C^{k}_{diff, 0}(M)$, we denote
the space of differential Hochschild $k$-cochains on $C^{\infty}(M)$
(i.e. k-multidifferential operators on $M$) vanishing on constants,
and denote by $b: C^{k}_{diff, 0}(M) \lon C^{k+1}_{diff, 0}(M) $
the Hochschild coboundary operator. 

\begin{pro} 
Suppose that $*_{\hbar}$ and $*'_{\hbar}$ are two compatible
star-products on $M$:
$$u*_{\hbar}v=\sum_{k=0}^{\infty} \hbar^{k} C_{k}(u, v), \ \ \ 
u*'_{\hbar}v=\sum_{k=0}^{\infty} \hbar^{k} C'_{k}(u, v), \ \forall u, v
\in C^{\infty}(M).$$
Assume that $*_{\hbar}$ and $*'_{\hbar}$ coincide with each other  up
to order $n$, i.e., $C_{k}  =C'_{k}, \ 0\leq k\leq n$. Then
\begin{enumerate}
\item $(C_{n+1}- C'_{n+1})(u, v)=\Vec{B}(u, v)+(b \Vec{E})(u, v)$,
where $B\in C^{2}(\frakh^* , (\wedge^{2}\frakg )^{\frakh})$ is
a $\delta_r$ 2-cocycle (i.e., $\delta_r B=0$), and  $E: \frakh^* 
\lon (U\frakg )^{\frakh}$. Here $\delta_r$  denotes the
coboundary operator defined by Equation (\ref{eq:deltar}). 
\item $C_{1}=\half \{\cdot, \cdot\}+b\Vec{c_{1}}$ for some $c_{1}\in 
C^{\infty} (\frakh^* , (U\frakg )^{\frakh} )$.
\item If $B=\delta_r  X, \ X\in  C^{\infty} (\frakh^* , \frakg^{\frakh})$,
then the formal operator $T=1+\hbar^{n}\Vec{X}
+\hbar^{n+1}\Vec{E_1}$ transforms $*_{\hbar}$  to another
star-product, which coincides with $*'_{\hbar}$ up to order
$n+1$. Here ${E_1}={E}(u)-[{X}, c_{1}]$. 
\end{enumerate}
\end{pro}
\pf We use a similar argument as in \cite{BBG}.

(i). By definition, if either $u$ or $v$ is in $C^{\infty}(\frakh^* )$, we have
$u*_{\hbar}v=u*'_{\hbar}v=\vT (u, v)$, which  implies that
$ (C_{n+1}- C'_{n+1})(u, v)=0$.

On the other hand, as it  is well known,
 $C_{n+1}- C'_{n+1}$ is a \Hoch 2-cocycle \cite{BBG, Xu0}. Hence we may write
$$C_{n+1}- C'_{n+1}=S+bT, $$
where $S\in \gm (\wedge^{2}TM)$ and $T$ is a \Hoch 1-cochain. Since
$S$ and $bT$ are, respectively,  the skew-symmetric  and
 symmetric parts of $C_{n+1}- C'_{n+1}$, they share 
many common properties as $C_{n+1}- C'_{n+1}$.
In particular, both of them  are $G\times H$-invariant and vanish
when one of the argument  $u$ or $v$ belongs to  $C^{\infty}(\frakh^* )$. 
This implies that $S=\Vec{B}$, for some $B\in C^{\infty}(\frakh^* ,
(\wedge^{2}\frakg )^{\frakh})$. It is also standard \cite{BBG, Xu0} 
that $S$ satisfies the equation:
$[\pi , S]=0$, which is equivalent to $\delta_r  B=0$ according
to the remark following  Proposition \ref{pro:delta}.

Now $M=\frakh^* \times G$ clearly 
admits a $G\times H$-invariant  (in fact  G-biinvariant)
connection. 
Since $bT$ is $G\times H$-invariant, according to Proposition 2.1 
in \cite{BBG}, we can  assume that $T$ is a $G\times H$-invariant
1-cochain.
Since $(bT)(u, v)=0$,  $\forall u, v\in C^{\infty}(\frakh^* )$,
we have   $u(Tv)-T(uv)+(Tu)v=0$.
On the other hand, since $Tu$ is $G$-invariant, it must be
a function of $\lambda \in \frakh^*$ only, i.e.,  $Tu \in C^{\infty}(\frakh^* )$.
Hence  the restriction of the operator $T$ to
$ C^{\infty}(\frakh^* )$  defines a vector field $Y $ 
on $\frakh^* $.
 Now since $(bT) (u, v)=0, \forall u\in C^{\infty}(\frakh^* )$,
it follows  that
$$(T-Y)(uv)=u(T-Y)(v), \ \ \ \ \forall u\in  C^{\infty}(\frakh^* ), \
v\in C^{\infty}(M).$$
Hence $T-Y$ does not involve any derivative with respect to
 $\lambda\in \frakh^*$.
Since $T-Y$ is $G\times H$-invariant,   it follows that
$T-Y=\Vec{E}$, for  some $E: \frakh^*
\lon (U\frakg)^{\frakh}$.
Therefore, $bT =b\Vec{E}$.

(ii) It is standard that  $C_{1}=\half \{\cdot, \cdot\}+ bc_{1}'$,
where $c_{1}'$ is a Hochschild 1-cochain.
By repeating a similar argument as in (i), we 
can prove that $c_{1}'$ can be chosen so that $c_{1}'=\Vec{c_{1}}$ 
for some $c_{1}\in
C^{\infty} (\frakh^* , (U\frakg )^{\frakh} )$.

(iii) If $B=\delta_r  X$, then $\Vec{B}=[\pi , \Vec{X}]$ according
 to  the remark following Proposition \ref{pro:delta}.
 It is easy to check that the operator
$T=1+\hbar^{n}\Vec{X}
+\hbar^{n+1}\Vec{E_1}$ transforms $*_{\hbar}$  to another
star-product, which coincides with $*'_{\hbar}$ up to order
$n+1$.
\qed

As a consequence, we have

\begin{cor}
\label{cor:eq-F}
If $r$ is a non-degenerate triangular  dynamical $r$-matrix
and $M=\frakh^* \times G$ its associated  symplectic manifold,
then every compatible $*$-product on $M$ is
strongly equivalent to a Fedosov $*$-product as constructed
in Corollary  \ref{cor:compatible}.
\end{cor}
\pf This follows essentially from the same
  argument as in the proof of Proposition 4.1 in \cite{BBG}.
 We will omit it here. \qed

Combing with   Theorem \ref{thm:Fed-classification}, we thus have
proved:

\begin{thm}
Let $M=\frakh^* \times G$ be the symplectic manifold
associated  with  a non-degenerate triangular  dynamical $r$-matrix
$r:\frakh^* \lon \wedge^2\frakg$. Then the equivalent classes of 
compatible $*$-products on $M$  are classified
by the relative Lie algebra cohomology (with coefficients
being formal power series of $\hbar$) $H^{2}(\frakg , \frakh )\flb\hbar \frb$. 
\end{thm}

Using Theorem \ref{thm:transformation}, we  are thus lead to
the following 

\begin{thm}
The equivalence classes of quantization of a non-degenerate triangular
 dynamical $r$-matrix 
$r:\frakh^* \lon \wedge^2  \frakg$ are classified by
the  relative Lie algebra cohomology (with coefficients
being formal power series of $\hbar$) $H^{2}(\frakg , \frakh )\flb\hbar \frb$.
\end{thm}
{\bf Remark} It would be interesting to see if  this 
theorem  can be proved by  directly applying
the usual classification theorem of star products on a symplectic manifold.
One of the difficulties is that the characteristic
class of a star product is usually difficult to computer.
 Recently, Tsygan comes up a
nice way of redefining the characteristic
class using the jet bundle. This may shed some new light on our problem.

Inspired by Kontesvich's formality theorem, we end
this section with the following:\\\\\\
{\bf Conjecture} For an arbitrary   classical triangular  dynamical $r$-matrix
$r:\frakh^* \lon \wedge^{2}\frakg$, the quantization is
classified by $\calm_{r} (\frakg \flb\hbar \frb, \frakh )$, the
formal neighbourhood  of $r$ in the moduli space $\calm 
(\frakg \flb\hbar \frb, \frakh )$. 

\appendix
\section{Appendix}

In this section, we  recall some basic ingredients  of
the Fedosov construction of
$*$-products on a symplectic manifold, as well as some
useful notations,   which are used throughout the paper.
For details, readers should
consult \cite{F1, F2}.

Let $(M, \omega )$ be a symplectic manifold of dimension $2n$.
Then, each tangent space $T_{x}M$ is equipped with  a linear symplectic
structure, which can be quantized using  the standard Moyal-Weyl product.
The resulting space is denoted by $W_{x}$.  More precisely,

\begin{defi}
A formal Weyl  algebra $W_{x}$  associated to $T_{x}M$
is an associative algebra  with a unit over $\complex $,
whose elements consist of formal  power  series in $\hbar $ with coefficients
being formal  polynomials in $T_{x}M$. In other words, each element has the form:
\begin{equation}
\label{eq:general}
a(y, \hbar )=\sum \hbar^{k}a_{k, \alpha }y^{\alpha }
\end{equation}
where   $y=(y^{1}, \cdots , y^{2n})$ is a linear coordinate system
on $T_{x}M$, $\alpha =(\alpha_{1}, \cdots , \alpha_{2n})$
is a multi-index,       $y^{\alpha}=(y^{1})^{\alpha_{1}}\cdots
(y^{2n})^{\alpha_{2n}}$, and $a_{k, \alpha }$ are constants.
The product is defined according to the Moyal-Weyl rule:

\begin{equation}
\label{eq:moyal}
a  *b=\sum_{k=0}^{\infty} (\frac{\hbar}{2})^{k}\frac{1}{k!}\pi^{i_{1}j_{1}}
\cdots \pi^{i_{k}j_{k}}
\frac{\partial^{k} a}{\partial y^{i_{1}}\cdots \partial y^{i_{k}}}
\frac{\partial^{k} b}{\partial y^{j_{1}}\cdots \partial y^{j_{k}}} .
\end{equation}
\end{defi}


Let $W=\cup_{x\in M}W_{x}$. Then $W$ is a bundle of algebras over $M$, called
the Weyl bundle.
 Its  space of   sections $\gm W$  forms  an associative algebra
with unit under the fiberwise multiplications.
One may think of $W$ as a ``quantum tangent bundle"
of  $M$, whose space of sections $\gm W$ gives
rise to a deformation quantization for the
tangent bundle $TM$, considered as a Poisson manifold
with fiberwise linear symplectic structures. As in \cite{F1},
 by $W^{+}$ we denote the extension of the algebra $W$ consisting of
those elements described  as follows:
\begin{enumerate}
\item elements  $a\in W^{+} $ are given by series (\ref{eq:general}), but
the powers of $\hbar$ can be both positive and negative;
\item the total degree $2k +|\alpha |$ of any term of the
series is nonnegative;
\item there exists a finite number of terms with a given
nonnegative total degree.
\end{enumerate}

 The center $Z (W)$  of $\gm W$ consists of sections
 not containing $y's$, thus can be naturally
identified with $C^{\infty}(M)\flb\hbar \frb$.
By assigning degrees to $y's$ and $\hbar$
with $\mbox{deg} y^{i}=1$ and $\mbox{deg} \hbar =2$, there
is a natural filtration
$$ C^{\infty}(M) 
\subset \gm (W_{1})\subset \cdots \gm (W_{i})\subset \gm (W_{i+1})\cdots
\subset \gm (W)    $$ with respect
to the total degree (e.g.,  any individual  term in  the summation of the RHS 
of Equation (\ref{eq:general}) has degree $2k+|\alpha|$.)

A differential $q$-form with values in $W$ is a section of the
bundle $W\otimes \wedge^{q}T^{*}M$, which  
can be expressed locally as
\begin{equation}
a(x, y, \hbar , dx)=\sum \hbar^{k}   a_{k, i_{1}\cdots i_{p}, j_{1}\cdots j_{q}}
y^{i_{1}}\cdots y^{i_{p}}dx^{j_{1}}\wedge \cdots \wedge dx^{j_{q}}.
\end{equation}
Here the coefficient $a_{k, i_{1}\cdots i_{p}, j_{1}\cdots j_{q}} $
is a covariant  tensor symmetric with respect to
$i_{1} \cdots i_{p}$ and antisymmetric in $j_{1}\cdots j_{q}$.
For short,  we denote the space of these
sections by  $\gm W\otimes \Lambda^{q}$. There is an associative
product $\smalcirc$ on $\gm W\otimes \Lambda^{*}$, which naturally  extends
the multiplication $*$ on $\gm W$ and the wedge product 
on $ \Lambda^{*}$:
\begin{equation}
\label{eq:circle}
(a\ot \theta )\smalcirc (b\ot \omega )=(a*b)\ot (\theta \wedge \omega ), 
\ \ \ \forall a, b\in \gm W , \ \mbox{ and } \  \theta , \omega \in \Lambda^{*}.
\end{equation}

The usual   exterior derivative on differential forms
 extends, in a straightforward way,
 to an operator $\delta$ on $W$-valued
differential forms:

\begin{equation}
\delta a=\sum_{i} dx^i\wedge \frac{\partial a}{\partial y^i}, \ \ \ \forall a\in \gm W\otimes \Lambda^* .
\end{equation}

By $\delta ^{-1}$, we denote its ``inverse" operator defined  by:
\begin{equation}
\delta ^{-1}a = \sum_{i} \frac{1}{p+q}y^{i}
 ( \frac{\partial }{\partial x^i}\per a)
\end{equation}
 when  $p+q>0$, and $\delta ^{-1}a=0$ when $p+q=0$,
where    $a\in \gm W\otimes \Lambda^{q}$ is
  homogeneous of degree $p$ in $y$.

There is a  ``Hodge"-decomposition:
\begin{equation}
\label{eq:hodge}
a=\delta \delta ^{-1}a+\delta ^{-1}\delta a+a_{0 0}, \ \ \ \forall a\in \gm W\otimes \Lambda^{*},
\end{equation}
where $a_{00}(x)$ is the  constant term of $a$, i.e, the $0$-form
term of $a|_{y=0}$ or  $a_{00}(x)=a(x, 0, 0, 0)$.
The operator $\delta$  possesses most of the 
basic properties of the usual exterior 
 derivatives.  For example,
$$\delta^{2}=0 \ \ \mbox{and } (\delta^{-1})^{2}=0 .$$

It is also clear that both $\delta$ and $\delta^{-1}$ commute
with the Lie derivative, i.e., $\forall X\in {\frak X}  (M)$,
\begin{eqnarray}
L_{X}\smalcirc \delta = \delta \smalcirc L_{X}, \ \ \ \mbox{ and }\ 
L_{X}\smalcirc  \delta^{-1}=\delta^{-1}\smalcirc  L_{X}.
\end{eqnarray}

Let $\nabla$ be a torsion-free symplectic connection on $M$ and 
$$\partial : \gm W \lon \gm W \otimes \Lambda^{1}$$
 be its induced covariant
derivative.

Consider   a connection on $W$ of the form:
\begin{equation}
\label{eq:connection}
D=-\delta +\partial + \ih [\gamma ,\cdot \ ],
\end{equation}
with $\gamma \in \gm W \otimes \Lambda^{1}$.

Clearly,  $D$ is a derivation with respect to the  Moyal-Weyl product, i.e.,
\begin{equation}
D(a* b)=a * Db +Da*  b, \  \ \ \forall a, b\in \gm W .
\end{equation}

 A simple calculation yields that
\begin{equation}
D^{2}a=-[\ih \Omega , a] , \ \forall a\in \gm W,
\end{equation}
where
\begin{equation}
\label{eq:curvature}
\Omega =\omega -R+\delta \gamma -\partial \gamma -\ih \gamma^{2}.
\end{equation}
Here  $R=\frac{1}{4}R_{ijkl}y^{i}y^{j}dx^{k}\wedge dx^{l}$,  and
$R_{ijkl}=\omega_{im}R^{m}_{jkl}$ is the curvature tensor
of the symplectic connection as defined by  Equation  (\ref{eq:sym-curvature}).

A connection  of the form (\ref{eq:connection}) 
 is called {\em Abelian} if
$\Omega $ is a scalar 2-form, i.e., $\Omega \in \Omega^{2}(M)\flb\hbar \frb$.
It  is called  a {\em Fedosov
 connection } if it is Abelian and in addition $\gamma \in \gm W_{3}\otimes \Lambda^{1}$.
For an Abelian connection, the Bianchi identity  implies that
$d\Omega =D\Omega =0$, i.e., $\Omega \in Z^{2}(M)\flb\hbar \frb$.
In this case,  $\Omega $ is called the Weyl curvature.

\begin{thm} (Fedosov \cite{F2})
\label{label:f1}
Let  $\nabla  $ be   a torsion-free symplectic connection,
 and $\Omega =\omega +\hbar \omega_{1} +\cdots \in
Z^{2}(M)\flb\hbar \frb$  a  perturbation of the symplectic form
in the space  $Z^{2}(M)\flb\hbar \frb$.
 There exists
a unique $\gamma \in \gm W_{3}\otimes \Lambda^{1}$
such that $D$,  given by Equation (\ref{eq:connection}),
 is a Fedosov connection, which has $\Omega$ as the Weyl curvature  and 
  satisfies
 $$\delta^{-1}\gamma =0.$$
\end{thm}
\pf It suffices to solve  the equation:
\begin{equation}
\omega -R+\delta \gamma -\partial \gamma -\ih \gamma^{2}=\Omega.
\end{equation}
This  is equivalent to
\begin{equation}
\label{eq:r}
\delta \gamma =\tilde{\Omega}+\partial \gamma +\ih \gamma^{2},
\end{equation}
where 
\begin{equation}
\tilde{\Omega}=\Omega -\omega+R.
\label{eq:tildeOme}
\end{equation}
Applying the operator $\delta^{-1}$ to Equation (\ref{eq:r})
and using the Hodge decomposition  (Equation (\ref{eq:hodge})), we obtain
\begin{equation}
\label{eq:r-iteration}
\gamma =\delta^{-1}\tilde{\Omega}+\delta^{-1}(\partial \gamma +\ih \gamma^{2} ).
\end{equation}
Note that $\gamma_{00}=0$ since $\gamma$ is a $1$-form.

Take  $\gamma_{0}=\delta^{-1}\tilde{\Omega}$, and consider the following
iteration equation: 
\begin{equation}
\label{eq:r-iteration2}
\gamma_{n+1} =\gamma_{0}+\delta^{-1}(\partial \gamma_{n} +\ih \gamma_{n}^{2} ),
\ \ \ \ \forall n\geq 0.
\end{equation}

Since the operator $\partial$ preserves the filtration 
and $\delta^{-1}$ raises it by
$1$, $\gamma_{n}$ defined by  Equation (\ref{eq:r-iteration2}) 
converges  to a unique $\gamma \in \gm W \otimes \Lambda^1$, which is clearly
a solution to Equation (\ref{eq:r-iteration}).
 Moreover
since $\gamma_{0}$ is at least  of degree 3, 
$\gamma$ is indeed an element   in  $\gm W_{3}\otimes \Lambda^{1}$. \qed
 
Theorem \ref{label:f1}  indicates that a Fedosov
connection $D$ is uniquely determined by a torsion-free symplectic
connection $\nabla$ and a Weyl curvature 
$\Omega=\sum_{i=0}^{\infty}\hbar^{i}\omega_{i}\in
Z^{2}(M)\flb\hbar \frb$. For this reason, we will say  that  $D$
is a Fedosov connection corresponding to the pair
$(\nabla , \Omega )$.

If $D$ is a Fedosov connection,  the space of all parallel
sections $W_{D}$ automatically becomes  an  associative
algebra. 
Let $\sigma$ denote the projection from $W_{D}$ to its center
 $C^{\infty}(M)\flb\hbar \frb$ defined by  $\sigma (a)=a|_{y=0}$.

\begin{thm} (Fedosov \cite{F2})
\label{thm:f2}
For any  $a_{0}(x, \hbar)\in C^{\infty}(M)\flb\hbar \frb$ there is a unique
section $a\in W_{D}$ such that $\sigma (a)=a_{0}$. Therefore,
$\sigma$ establishes an isomorphism between $W_{D}$ and
 $C^{\infty}(M)\flb\hbar \frb$ as vector spaces.
\end{thm}
\pf 
The equation $Da=0$  can be written as
$$\delta a=\partial a +[\ih \gamma , a].$$
Applying the operator $\delta^{-1}$, it follows from the Hodge
decomposition  (Equation (\ref{eq:hodge}))  that
\begin{equation}
\label{eq:parallel}
a=a_{0}+\delta^{-1}(\partial a+[\ih \gamma , a]).
\end{equation}
In analogue to the proof of Theorem \ref{label:f1},
we can solve this equation 
 by the iteration formula:
\begin{equation}
a_{n+1} =a+\delta^{-1}(\partial {a_{n}} +[\ih \gamma , {a_{n}} ]).
  \end{equation} \qed

\end{document}